\documentclass[12pt,draftcls,onecolumn]{IEEEtran}

\usepackage{etoolbox}
\newtoggle{draft}
\toggletrue{draft}
\iftoggle{draft}{%

}{%

}

\IEEEoverridecommandlockouts                              
\usepackage{graphics,cite} 
\usepackage{epsfig,epstopdf} 
\usepackage{amsmath} 
\usepackage{amssymb}  
\usepackage{url}
\usepackage{hyperref}
\usepackage{enumerate}
\usepackage{subfigure}
\usepackage{algorithmic}
\usepackage{algorithm}
\usepackage{enumitem}
\usepackage{array}

\DeclareGraphicsExtensions{.eps,.pdf,.png}

\title{Optimal scaling of the ADMM algorithm for\\ distributed quadratic programming}

\author{Andr\'{e} Teixeira, Euhanna Ghadimi, Iman Shames,\\ Henrik Sandberg, and Mikael Johansson 
\thanks{A.~Teixeira, E.~Ghadimi, H.~Sandberg, and M.~Johansson are with the ACCESS Linnaeus Centre, Electrical Engineering, KTH Royal Institute of Technology, Stockholm, Sweden.
{\tt\small \{andretei,euhanna,hsan,mikaelj\}@kth.se}}
 \thanks{I.~Shames is with the Department of Electrical and Electronic Engineering, University of Melbourne, Australia. ~{\tt\small iman.shames@unimelb.edu.au.}}
 \thanks{This work was sponsored in part by the Swedish Foundation for Strategic Research, the Swedish Research Council and a McKenzie Fellowship.}
}

\begin{document}
\newtheorem{proposition}{Proposition}
\newtheorem{lem}{Lemma}
\newtheorem{theorem}{Theorem}
\newtheorem{thm}{Theorem}
\newtheorem{rem}{Remark}
\newtheorem{assumption}{Assumption}
\newtheorem{example}{Example}
\newtheorem{definition}{Definition}
\newtheorem{problem}{Problem}
\newtheorem{cor}{Corollary}
\newcommand{\ie}{i.e.}
\newcommand{\eg}{e.g.}
\newcommand{\cf}{cf.}
\newcommand{\R}[1]{\mathbf{R}^{#1}}
\newcommand{\PD}[1]{\mathcal{S}_{++}^{#1}}
\newcommand{\PSD}[1]{\mathcal{S}_{+}^{#1}}
\newcommand{\Null}[1]{\mathcal{N}(#1)}
\newcommand{\Range}[1]{\mathcal{R}(#1)}
\newcommand{\PrjR}[1]{\Pi_{\Range{#1}}}
\newcommand{\PrjN}[1]{\Pi_{\Null{#1}}}

\newcommand{\MJ}[1]{{\color{black}{#1}}}
\newcommand{\EG}[1]{{\color{black}{#1}}}
\newcommand{\AT}[1]{{\color{black}{#1}}}
\newcommand{\IS}[1]{{\color{black}{#1}}}


\maketitle
\begin{abstract}                          
This paper presents optimal scaling of the alternating directions method of multipliers (ADMM) algorithm for a class of distributed quadratic programming problems. The scaling corresponds to the ADMM step-size and relaxation parameter, as well as the edge-weights of the underlying communication graph. We optimize these parameters to yield the smallest convergence factor of the algorithm. Explicit expressions are derived for the  step-size and relaxation parameter, as well as for the corresponding convergence factor. Numerical simulations justify our results and highlight the benefits of optimally scaling the ADMM algorithm.
\end{abstract}
\section{Introduction}
\label{sec:introduction}

Recently, a number of applications have triggered a strong interest in distributed algorithms for large-scale quadratic programming. These applications include multi-agent systems~\cite{nedic10,Italian}, distributed model predictive control~\cite{accelerated_mpc13, farhad12}, and state estimation in networks~\cite{Falcao1995}, to name a few.
As these systems become larger and their complexity increases, more efficient algorithms are required. It has been argued that
the alternating direction method of multipliers (ADMM) is a particularly powerful  approach~\cite{Boyd11}.
One attractive feature of ADMM is that it is guaranteed to converge for all (positive) values of its step-size parameter. This contrasts many alternative techniques, such as dual decomposition, where mistuning of the step-size for the gradient updates can render the iterations unstable.

The ADMM method has been observed to converge fast in many applications~\cite{Boyd11,tyler12,marriette12, Mota12} and for certain classes of problems it is known to converge at a linear rate~\cite{luo12,boley12, deng12}. However, the solution times are sensitive to the choice of the step-size parameter, and when this parameter is not properly tuned, the ADMM iterations may converge (much) slower than the standard gradient algorithm~\cite{GTS:13}.
In practice, the ADMM algorithm parameters are tuned empirically for each specific application. For example,~\cite{tyler12,marriette12, Mota12} propose different rules of thumb for picking the step-size for different distributed quadratic programming applications, and empirical results for choosing the best relaxation parameter can be found in~\cite{Boyd11}. However, a thorough analysis and design of the \emph{optimal} step-size, relaxation parameter, and scaling rules for the ADMM algorithm is still missing in the literature.

The aim of this paper is to address this shortcoming by deriving jointly optimal ADMM parameters for a class of distributed quadratic programming problems that appears in applications such as distributed power network state-estimation~\cite{Gomez2011_DSE} and distributed averaging~\cite{Italian}. In this class of problems, a number of agents collaborate with neighbors in a graph to minimize a convex objective function over a combination of shared and private variables.
By introducing local copies of the global decision vector at each node, unconstrained quadratic programming problems in this class can be re-written as equality-constrained quadratic programming problems, where the constraints enforce consistency among the local decision vectors. By analyzing these equality-constrained quadratic programming problems, we are able to characterize the optimal step-size, over-relaxation and constraint scalings for the associated ADMM iterations.

Specifically, since the ADMM iterations for our problems are linear, the convergence behavior depends on the spectrum of the transition matrix. In each step, the distance to the optimal point is guaranteed to decay by a factor equal to the largest non-unity magnitude eigenvalue. We refer to this quantity as the {\it convergence factor}.
We show that the eigenvalues of the transition matrix are given by the roots of quadratic polynomials whose coefficients depend on the step-size, relaxation parameter, and the spectrum of the graph describing interactions between agents. Several properties of the roots with respect to the ADMM parameters are analyzed and used to develop scaled ADMM iterations with a minimal convergence factor. Analytical expressions for the proposed step-size, relaxation parameter, and the resulting convergence factor are derived. Finally, given that the optimal step-size and relaxation parameter are chosen, we propose methods to further improve the convergence factor by optimal scaling. The optimal step-size for the standard ADMM iterations (without the relaxation parameter) was characterized in prior related work~\cite{teixeira_admm_2013}, while a brief summary of the results and their application to distributed averaging problems are reported in~\cite{ghadimi2014admm}.

The outline of this paper is as follows. Section~\ref{sec:background} gives a background on the ADMM and illustrates how the ADMM may be used to formulate distributed optimization problems as equality-constrained optimization problems. The ADMM iterations for equality-constrained quadratic programming problems are formulated and analyzed in Section~\ref{sec:qp_equality}. Distributed quadratic programming and optimal networked-constrained scaling of the ADMM algorithm are addressed in Section~\ref{sec:distributed_QP}. Numerical examples illustrating our results and comparing them to state-of-the art techniques are presented in Section~\ref{sec:numerical}. Section~\ref{sec:conclusion} concludes the paper.


\section{Background}
\label{sec:background}

Below we define the notation used throughout the paper, followed by a summary of the ADMM algorithm and its application to distributed optimization problems.

\subsection{Notation}
The cardinality of a set $\mathcal{A}$ is expressed as $\vert \mathcal{A}\vert$.  The sets of real and complex numbers are denoted by $\mathbf{R}$ and $\mathbf{C}$, respectively. The dimension of a subspace $\mathcal{X}$ is denoted by $\mbox{dim}(\mathcal{X})$ and for a given matrix $A$, $\mbox{span}(A)$ is the subspace spanned by its columns.  For $A\in \mathbf{R}^{n\times m}$, $\Range{A}\triangleq \{y\in\R{n} \vert \; y=Ax,\, x\in\mathbf{R}^{m}  \}$ denotes its range-space and $\Null{A}\triangleq\{x\in\mathbf{R}^m \vert \; Ax=0\}$ its null-space. For $A$ with full-column rank, $A^\dagger \triangleq (A^\top A)^{-1}A^\top$ is the pseudo-inverse of $A$ and $\Pi_{\Range{A}} \triangleq A A^\dagger$ is the orthogonal projector onto $\Range{A}$.
Consider $B,D\in \R{n \times n}$, with $D$ being invertible. The generalized eigenvalues of $(B, D)$ are defined as the values $\lambda \in \mathbf{C}$ such that  $(B - \lambda D) v = 0$ holds for some nonzero vector $v\in\mathbf{C}^n$. The set of real-symmetric matrices in $\mathbf{R}^{n\times n}$ is denoted by $\mathcal{S}^n$. Additionally, $A\succ 0$ ($A\succeq 0$) indicates that $A$ is positive definite (semi-definite). Given a sequence of $m$ square matrices $\left\{ A_i \right\}_{i=1}^m$ with $A_i\in\mathbf{R}^{n\times n}$, we denote $\mbox{diag}\left( \left\{ A_i \right\}_{i=1}^m \right)\in\mathbf{R}^{nm \times nm}$ as the block-diagonal matrix with $A_i$ in its $i$-th diagonal block.
Given a vector $x\in\mathbf{R}^n$, its Euclidean norm is denoted as $\| x \|_2 = \sqrt{x^\top x}$. The convergence factor of a sequence of vectors $\{ \sigma^k\}$, with $\sigma^k \in\mathbf{R}^n$ for all $k$, converging to $\sigma^{\star}\in\mathbf{R}^{n}$ is defined as
\begin{equation}\label{eq:convergence_factor}
\begin{aligned}
\phi^\star &\triangleq
\limsup_{k\rightarrow \infty}
\dfrac{\|\sigma^{k+1} - \sigma^\star\|_2}{\|\sigma^{k} - \sigma^\star\|_2}.
\end{aligned}
\end{equation}

Let ${\mathcal{G}}(\mathcal{V},\mathcal{E})$ be an undirected graph with vertex set $\mathcal{V}$  and edge set $\mathcal{E}$.
For any given ordering of the edges of ${\mathcal{G}}$, the $k$-th edge is denoted by $e_k\in \mathcal{E}$.
Let ${\mathcal N}_i\triangleq\{j\neq i \vert \{i,j\}\in \mathcal{E}\}$ be the neighbor set of node $i$. Moreover, the sparsity pattern induced by $\mathcal{G}$ is defined as $\mathcal{A} \triangleq \{S \in \mathcal{S}^{\vert \mathcal{V}\vert} \vert S_{ij}=0 \,\mbox{if} \, i \neq j \, \mbox{and}\, \{i,j\}\,\not\in \mathcal{E} \}$. Given the undirected graph $\mathcal{G}(\mathcal{V},\mathcal{E})$, we introduce an associated weighted directed graph $\bar{\mathcal{G}}(\mathcal{V},\bar{\mathcal{E}}, \mathcal{W})$. The edge set $\bar{\mathcal{E}}$ of $\bar{\mathcal{G}}$ contains two directed edges $(i,j)$ and $(j,i)$ for each undirected edge $\{i,j\} \in \mathcal{E}$, and the edge weights $\mathcal{W}=\{ W_{(i,j)} \}_{(i,j)\in \bar{\mathcal{E}}}$ comprise matrix-valued weights $W_{(i,j)}\in\mathbf{R}^{n_w\times n_w}$ with $W_{(i,j)}\succeq 0$ for each directed edge $(i,j)\in \bar{\mathcal E}$. Similarly to $\mathcal{G}$, we assume there exists an arbitrary ordering of the edges of $\bar{\mathcal{G}}$ in which $\bar{e}_k\in \bar{\mathcal{E}}$ denotes the $k$-th edge. The matrix $B^+\in\mathbf{R}^{|\bar{\mathcal{E}}|\times \vert \mathcal{V}\vert}$ is defined as $B^+_{kj}=1$ if $j$ is the head of $\bar{e}_k\in\bar{\mathcal{E}}$ and $B^+_{kj}=0$ otherwise. Similarly, define $B^-\in\mathbf{R}^{|\bar{\mathcal{E}}|\times  \vert \mathcal{V}\vert}$ so that $B^-_{kj}=1$ if $j$ is the tail of $\bar{e}_k\in\bar{\mathcal{E}}$ and $B^-_{kj}=0$ otherwise. Moreover, the edge-weight matrix is defined as $W\triangleq\mbox{diag}\left(\left\{W_{\bar{e}_k}\right\}_{i=1}^m \right)$. In the following sections, $\bar{\mathcal{G}}$ is used to describe scenarios where nodes $i$ and $j$ assign different weights to the undirected edge $\{i,j\}$.

For $n_w =1$ and symmetric weights $W_{\{i,j\}} \triangleq W_{(i,j)} = W_{(j,i)}>0$, $\bar{\mathcal{G}}$ is equivalent to a weighted undirected graph whose adjacency matrix $A\in\mathcal{A}$ is defined as $A_{ij} = W_{\{i,j\}}$ for $\{i,j\}\in\mathcal{E}$ and $A_{ii}=0$.
The corresponding diagonal degree matrix $D$ is given by $D_{ii} = \sum_{j\in\mathcal{N}_i}A_{ij}$.

\subsection{The ADMM method}
The ADMM algorithm solves problems of the form
\begin{align}\label{eq:constrained problem}
	\begin{array}[c]{ll}
	\underset{x,\,z}{\mbox{minimize}} & f(x)+g(z)\\
	\mbox{subject to} & Ex+Fz - h = 0
	\end{array}
\end{align}
where $f$ and $g$ are convex functions, $x\in {\mathbf R}^n$, $z\in {\mathbf R}^m$, $h\in {\mathbf R}^p$. Moreover, $E\in {\mathbf R}^{p\times n}$ and $F\in {\mathbf R}^{p\times m}$ are assumed to have full-column rank; see~\cite{Boyd11} for a detailed review.
The method is based on the \emph{augmented Lagrangian}
\begin{equation}
\begin{aligned}
\label{eq:augmented_Lagrangian}
\iftoggle{draft}{%
L_{\rho}(x,z,\mu) = f(x)+g(z) + (\rho/2)\Vert Ex+Fz - h \Vert_2^2 + \mu^{\top}(Ex+Fz - h)
}{%
L_{\rho}(x,z,\mu) =& f(x)+g(z) + (\rho/2)\Vert Ex+Fz - h \Vert_2^2  \\
                             &+\mu^{\top}(Ex+Fz - h)
}
%
%
\end{aligned}
\end{equation}
and performs sequential minimization of the $x$ and $z$ variables, followed by a dual variable update.
It is convenient to use the scaled dual variable $u=\mu/\rho$, which yields the iterations
\begin{align}
	x^{k+1} &= \underset{x}{\operatorname{argmin}}\; f(x)+ (\rho/2)\Vert Ex+Fz^k -h + u^k\Vert_2^2\nonumber \\
	z^{k+1} &= \underset{z}{\operatorname{argmin}}\; g(z) + (\rho/2)\Vert Ex^{k+1}+Fz - h +u^{k}\Vert_2^2  \label{eqn:admm_scaled}\\
	u^{k+1} &= u^{k} + Ex^{k+1} + Fz^{k+1} - h\nonumber
\end{align}

These iterations indicate that the method is particularly useful when the $x$- and $z$-minimizations can be carried out efficiently (\eg,~when they admit closed-form expressions). One advantage of the method is that there is only one single algorithm parameter, $\rho$, and that under rather mild conditions, the method can be shown to converge for all values of this parameter; see, \eg,~\cite{Boyd11}.
However, $\rho$ has a direct impact on the convergence speed of the algorithm, and inadequate tuning of this parameter may render the method very slow.

The convergence properties of iterative algorithms can often be improved by accounting for the past iterates when computing the next. This technique is called \emph{relaxation}. For ADMM it amounts to replacing $E x^{k+1}$ with $\gamma^{k+1} = \alpha^k E x^{k+1}- (1-\alpha^k) (Fz^k -h)$ in the $z$- and $u$-updates \cite{Boyd11}, yielding
\begin{equation}
\label{eqn:admm_relaxed}
\begin{aligned}
%
	z^{k+1} &= \underset{z}{\operatorname{argmin}} \; g(z) + \dfrac{\rho}{2}\left\Vert \gamma^{k+1} +Fz-h+u^{k}\right\Vert_2^2,  \\
	u^{k+1} &= u^{k} + \gamma^{k+1}+ Fz^{k+1}-h.
\end{aligned}
\end{equation}
The parameter $\alpha^k\in (0,2)$ is called the \emph{relaxation parameter}. Note that letting $\alpha^k=1$ for all $k$ recovers the original ADMM iterations~\eqref{eqn:admm_scaled}. Empirical studies have suggested that $\alpha^k>1$ (referred to as over-relaxation) is often advantageous and the guideline $\alpha^k\in [1.5, 1.8]$ has been proposed~\cite{Boyd11}.

In the remaining parts of this paper, we derive explicit expressions for the step-size $\rho$ and relaxation parameter $\alpha$ that minimize the convergence factor~\eqref{eq:convergence_factor} for a class of distributed quadratic programming problems. In terms of the standard form (\ref{eq:constrained problem}),  $g(z)$ is linear and $f(x)$ is quadratic with a Hessian matrix $Q\succ 0$ such that $Q=\kappa E^\top E$ for some $\kappa>0$. Table~\ref{tab:summary} summarizes the proposed choice of parameters.
Note that the parameters and the resulting convergence factor only depend on $\lambda_1$ and $\lambda_{n-s}$, where $s=\mbox{dim}\left(\mathcal{N}([ E \; F ])\right)$ and $\{\lambda_i\}_{i=1}^n$ are the generalized eigenvalues of $\left(E^\top (2\PrjR{F} - I ) E , E^\top E\right)$, ordered in increasing magnitude.

\renewcommand{\arraystretch}{1.5}
\begin{table}[H]
\caption{Optimized ADMM parameters for distributed QP (Theorem~\ref{thm:optimal_f_alpha}).
\label{tab:summary}}
\centering
\begin{tabular}{|c|cc|}
  \hline
   \textbf{Case}                                                               & ${\rho^\star }$                        & ${\alpha^\star}$ \\[0.1cm] \hline
   $\lambda_{n-s} \geq  |\lambda_1|$     & $\frac{1}{\sqrt{1-\lambda_{n-s}^2}}$   & $2$                 \\ 
  $|\lambda_1| > \lambda_{n-s} >0$
 & $\frac{1}{\sqrt{1-\lambda_{n-s}^2}}$   & $\frac{4\rho^\star+4}{2-\rho^\star(\lambda_{n-s} + \lambda_1 - 2 -\sqrt{\lambda_1^2 - \lambda_{n-s}^2}) }$     \\ 
  $0\geq \lambda_{n-s} \geq  \lambda_1$ & $1$                           & $\frac{4}{2-\lambda_1}$   \\[0.1cm]
  \hline
\end{tabular}
\end{table}
\renewcommand{\arraystretch}{1}

We highlight that the results in Table~\ref{tab:summary} may be directly applied to the case where the equality constraints in~\eqref{eq:constrained problem} are scaled by a matrix $R\in\mathbf{R}^{r\times p}$, yielding $R(Ex+Fz-h)=0$.
As shown in Section~\ref{sec:Optimal_const_scaling}, such scaling may be used to further improve of the convergence factor.

Next we describe a distributed unconstrained optimization problem that belongs to the class of problems considered in the paper and will be used as a motivating example in Section~\ref{sec:numerical}.

\subsection{ADMM for distributed optimization}\label{subsec:ADMM_DQP}
Consider a network of agents, each endowed with a local convex loss function $f_i(x)$, that collaborate to find the decision vector $x$ that results in the minimal total loss, \emph{i.e.}
\begin{align*}
	\begin{array}[c]{ll}
	\underset{x\in\mathbf{R}^{n_x}}{\mbox{minimize}} & \sum_{i\in\mathcal{V}} f_i(x).
	\end{array}
\end{align*}
The interactions among agents are described by an undirected graph $\mathcal{G}(\mathcal{V},\mathcal{E})$: agents are only allowed to share their current iterate with neighbors $j\in {\mathcal N}_i$ in ${\mathcal G}$. By introducing local copies $x_i\in\mathbf{R}^{n_x}$ of the global decision vector at each node $i\in\mathcal{V}$, the original problem can be re-written as an equality constrained optimization problem with decision variables $\{x_i\}_{i\in\mathcal{V}}$ and separable objective:
\begin{align}\label{eq:DQP_consensus}
	\begin{array}[c]{ll}
	\underset{\{x_{i}\}}{\mbox{minimize}} & \sum_{i\in \mathcal{V}} f_i(x_{i}) \\
\mbox{subject to} & x_{i} = x_{j},\quad \forall\,i, j \in \mathcal{V}.
	\end{array}
\end{align}
The equality constraints ensure that the local decision vectors $x_i$ of all agents agree at optimum.
Since the problem must be solved distributedly, we make the following assumption~\cite{NEO:09}.
\begin{assumption}
\label{assump:Graph_Connected}
The graph ${\mathcal{G}}(\mathcal{V},\mathcal{E})$ is connected.
\end{assumption}
When the communication graph is connected, all equality constraints in (\ref{eq:DQP_consensus}) that do not correspond to neighboring nodes in ${\mathcal G}$ can be removed without altering the optimal solution. The remaining inequality constraints can be accounted for in different ways as described next (cf.~\cite{Italian}).

\subsubsection{Enforcing agreement with edge variables}\label{subsec:Edge_variable}

One way to ensure agreement between the nodes is to enforce all pairs of nodes connected by an edge to have the same value,~\ie,~$x_i = x_j$ for all $\{i,j\}\in\mathcal{E}$. To include this constraint in the ADMM formulation, one can introduce
an auxiliary variable $z_{\{i,j\}}$ for each edge $\{i,j\}\in\mathcal{E}$. The local constraints $x_i=z_{\{i,j\}}$ and $x_j=z_{\{i,j\}}$ are then introduced for neighboring nodes $i$ and $j$, and an equivalent form of~\eqref{eq:DQP_consensus} is formulated as
\begin{equation}\label{eq:DQP_high}
\begin{aligned}
	\begin{array}[c]{ll}
	\underset{\{x_i\}, \{z_{\{i,j\}}\}}{\mbox{minimize}} &\sum_{i\in \mathcal{V}} f_i(x_i)\\
	\mbox{subject to} & R_{(i,j)}x_i = R_{(i,j)}z_{\{i,j\}},\quad\forall i\in \mathcal{V},\; \forall (i,j)\in\bar{\mathcal{E}}.
	\end{array}
\end{aligned}
\end{equation}
Here, $R_{(i,j)}\in\mathbf{R}^{n_x\times n_x}$ acts as a scaling factor for the constraint defined along each edge $ (i,j)\in\bar{\mathcal{E}}$ and $W_{(i,j)} \triangleq R_{(i,j)}^\top R_{(i,j)}\succeq 0$ is the weight of the edge $(i,j)\in\bar{\mathcal{E}}$. The edge weights $W_{(i,j)} $ are included to increase the degrees of freedom available for nodes to improve the performance of the algorithm. We will discuss optimal design of these constraint scalings in Section~\ref{sec:Optimal_const_scaling}. Note that when the edge variables are fixed, (\ref{eq:DQP_high}) is separable and each agent $i$ can find the optimal $x_i$ without interacting with the other agents.

The optimization problem~\eqref{eq:DQP_high} can be written in the ADMM standard form~\eqref{eq:constrained problem} as follows. Define $x=[x_1^\top\,\cdots\, x_{\vert \mathcal{V}\vert}^\top]^\top$, $z = [z_{e_1}^\top\, \cdots \, z_{e_{|\mathcal{E}|}}^\top]^\top$, $f(x) = \sum_{i\in \mathcal{V}} f_i(x_i)$ and recall the matrix $B^+$ defined in Section~\ref{sec:introduction}. Problem (\ref{eq:DQP_high}) can then be rewritten as
\begin{equation}\label{eq:consensus_problem_edge}
	\begin{array}[c]{ll}
	\underset{x,z}{\mbox{minimize}} &f(x)\\  
	\mbox{subject to} &
    REx + RFz = 0,
%
	\end{array}
\end{equation}
where
\begin{equation}\label{eq:consensus_problem_edge_matrices}
E =
B^+ \otimes I_{n_x},\quad
F = - \begin{bmatrix}
I_{|\mathcal{E}|}\\
I_{|\mathcal{E}|}
\end{bmatrix} \otimes I_{n_x},\quad
R =  \mbox{diag}(\{R_{\bar{e}_i}\}_{\bar{e}_i \in \bar{\mathcal{E}}}).
\end{equation}

\subsubsection{Enforcing agreement with node variables}\label{subsec:Node_variable}
Another way of enforcing the agreement among the decision makers is via node variables. In this setup, each agent $i$ has to agree with all the neighboring agents, including itself. In the ADMM formulation, this constraint is formulated as $x_i=z_j$ for all $j\in \mathcal{N}_i \cup \{ i\}$, where $z_i\in\mathbf{R}^{n_x}$ is an auxiliary variable created per each node $i$. The optimization problem can be written as
\begin{equation}\label{eq:DQP_node}
\begin{aligned}
  \begin{array}[c]{ll}
	\underset{\{x_i\},\{ z_i\}}{\mbox{minimize}} &\sum_{i\in \mathcal{V}} f_i(x_i)\\
	\mbox{subject to} & R_{(i,j)} x_i =  R_{(i,j)}z_j,\quad\forall i\in \mathcal{V},\; \forall j\in\{\mathcal{N}_i \cup \{i\}\},
	\end{array}
\end{aligned}
\end{equation}
where $R_{(i,j)}\in\mathbf{R}^{n_x\times n_x}$ and $W_{(i,j)} = R_{i,(i,j)}^\top R_{i,(i,j)}\succeq 0$  is the weight of the edge $(i,j)\in\bar{\mathcal{E}}$. Additionally, we also have $R_{(i,i)}\in\mathbf{R}^{n_x \times n_x}$ and define $W_{(i,i)} \triangleq R_{(i,i)}^\top R_{(i,i)}\succeq 0$ as the matrix-valued weight of the self-loop $(i,i)$.
Similarly to the previous section, recalling the matrices $B^+$ and $B^-$ defined in Section~\ref{sec:background}, the distributed quadratic problem~\eqref{eq:DQP_node} can be rewritten as~\eqref{eq:consensus_problem_edge}
with
\begin{equation}\label{eq:consensus_problem_node_matrices}
\iftoggle{draft}{%
E =  \begin{bmatrix}
B^+ \\
I_{\vert\mathcal{V}\vert}
\end{bmatrix}\otimes I_{n_x},\quad
F = - \begin{bmatrix}
B^- \\
I_{\vert\mathcal{V}\vert}
\end{bmatrix}\otimes I_{n_x} ,\quad
R =  \mbox{diag}( \{R_{\bar{e}_k}\}_{\bar{e}_k\in\bar{\mathcal{E}}} , \;  \{R_{(i,i)}\}_{i \in\mathcal{V}}).
}{%
  \begin{aligned}
  E &=  \begin{bmatrix}
B^+ \\
I_{\vert\mathcal{V}\vert}
\end{bmatrix}\otimes I_{n_x},\quad
F = - \begin{bmatrix}
B^- \\
I_{\vert\mathcal{V}\vert}
\end{bmatrix}\otimes I_{n_x} ,\\
R &=  \mbox{diag}( \{R_{\bar{e}_k}\}_{\bar{e}_k\in\bar{\mathcal{E}}} , \;  \{R_{(i,i)}\}_{i \in\mathcal{V}}).
  \end{aligned}
}
\end{equation}

\section{ADMM for equality-constrained quadratic programming problems}\label{sec:qp_equality}
In this section, we analyze and optimize scaled ADMM iterations for the following class of equality-constrained quadratic programming problems
\begin{equation}\label{eq:QP_problem}
	\begin{array}[c]{ll}
	\underset{x,z}{\mbox{minimize}} &\dfrac{1}{2}x^\top Qx + q^{\top}x+c^\top z \\  
	\mbox{subject to} &
    REx + RFz = R h.
	\end{array}
\end{equation}
where $Q\in \mathbf{R}^{n\times n}, Q\succ 0$, and $q\in \mathbf{R}^n$. We assume that $E$, and $F$ have full-column rank.
An important difference compared to the standard ADMM iterations described in the previous section is that the original constraints $Ex+Fz=h$ have been scaled by a matrix $R\in \mathbf{R}^{r\times p}$.
\begin{assumption}\label{assum:feasible}
The scaling matrix $R$ is chosen so that no non-zero vector $v$ of the form $v=Ex+Fz-h$ belongs to the null-space of $R$.
\end{assumption}
In other words, after the scaling with $R$, the feasible set in~\eqref{eq:constrained problem} remains unchanged. Letting $\bar{E}=RE$, $\bar{F}=RF$, and $\bar{h}=Rh$, 
the penalty term in the augmented Lagrangian becomes $\rho/2 \Vert \bar{E}x+\bar{F}z-\bar{h}\Vert^2$.

 Our aim is to find the optimal scaling that minimizes the convergence factor of the corresponding ADMM iterations. In the next lemma we show that \eqref{eq:QP_problem} can be cast to the more suitable form:
\begin{equation}\label{eq:QP_problem_suitable}
	\begin{array}[c]{ll}
	\underset{x,z}{\mbox{minimize}} &\dfrac{1}{2}x^\top Qx + p^{\top}x+c^\top z \\  
	\mbox{subject to} &
    REx + RFz = 0.
	\end{array}
\end{equation}
\begin{lem}\label{lem:QP_change_of_variable}
Let $(\hat{x}, \hat{z})$ and $(x^\star,z^\star)$ be any feasible solution and optimal solution to~\eqref{eq:QP_problem}, respectively. Then the optimization problem~\eqref{eq:QP_problem_suitable} has the optimal solution $(x^\star-\hat{x},z^\star-\hat{z})$ if the parameters $q$ and $p$ in~\eqref{eq:QP_problem} and~\eqref{eq:QP_problem_suitable} satisfy $p= Q\hat{x}+q$.
  \begin{IEEEproof}
  See Appendix~\ref{app:lem:QP_change_of_variable}.
  \end{IEEEproof}
\end{lem}
Without loss of generality we thus assume $\bar{h}=0$  in the remainder of the paper. The scaled ADMM iterations for~\eqref{eq:QP_problem} with fixed relaxation parameter $\alpha^k=\alpha$ for all $k$ then read 
\begin{equation}\label{eqn:admm_scaled_iterations}
\begin{aligned}
	x^{k+1} =& (Q+\rho \bar{E}^\top \bar{E})^{-1}\left(-q - \rho\bar{E}^\top(\bar{F}z^k + u^k) \right)  \\
	\iftoggle{draft}{%
z^{k+1} =& -(\bar{F}^\top \bar{F})^{-1}\left(\bar{F}^\top \left( \alpha\bar{E}x^{k+1} - (1-\alpha)\bar{F}z^k  + u^k \right) +c/\rho\right) \\
}{%
  z^{k+1} =& -(\bar{F}^\top \bar{F})^{-1}\bar{F}^\top \left( \alpha\bar{E}x^{k+1} - (1-\alpha)\bar{F}z^k  + u^k \right) \\
  &-(\bar{F}^\top \bar{F})^{-1} c/\rho\\
}
	u^{k+1} =& u^{k} + \alpha\bar{E}x^{k+1} - (1-\alpha)\bar{F}z^k + \bar{F}z^{k+1} .
\end{aligned}
\end{equation}
Inserting the expression for $z^{k+1}$ in the $u$-update yields
\begin{equation}\nonumber
u^{k+1}=\Pi_{\mathcal{N}(\bar{F}^\top)}\left( \alpha \bar{E}x^{k+1}  + u^k \right) -{\bar{F}(\bar{F}^\top\bar{F})^{-1}c}/{\rho}.
\end{equation}
Since $\mathcal{N}(\bar{F}^\top)$ and $\mathcal{R}(\bar{F})$ are orthogonal complements, this implies that $\Pi_{\mathcal{R}(\bar{F})} u^k = -{\bar{F}(\bar{F}^\top\bar{F})^{-1}c}/{\rho}$ for all $k$. Thus
\begin{equation}\label{eq:z_iterations}
\bar{F}z^{k+1} =  (1-\alpha)\bar{F}z^{k}  - \alpha\Pi_{\mathcal{R}(\bar{F})}\bar{E}x^{k+1}.
\end{equation}
By inserting this expression in the $u$-update and applying the simplified iteration recursively, we find that
\begin{equation}\label{eq:u_iterations}
u^{k+1}=
\Pi_{\mathcal{N}(\bar{F}^\top)}\left( u^0 + \alpha\sum_{i=1}^{k+1}\bar{E}x^{i}   \right)-{\bar{F}(\bar{F}^\top\bar{F})^{-1}c}/{\rho} .
\end{equation}
We now apply \eqref{eq:z_iterations} and \eqref{eq:u_iterations} to eliminate $u$ from the $x$-updates:
\begin{equation}
\label{eqn:x_k+1_0}
\begin{aligned}
x^{k+1} &=
  \alpha\rho(Q+\rho \bar{E}^\top \bar{E})^{-1}\bar{E}^\top\left(\Pi_{\mathcal{R}(\bar{F}^\top)}-\Pi_{\mathcal{N}(\bar{F}^\top)}\right)\bar{E}x^k\\
 &+ x^k + \alpha \rho(Q+\rho \bar{E}^\top \bar{E})^{-1}\bar{E}^\top\bar{F}z^{k-1}.
\end{aligned}
\end{equation}
Thus, using~\eqref{eqn:x_k+1_0}  and defining
$y^k\triangleq \bar{E}^\top \bar{F}z^k$, the ADMM iterations can be rewritten in the following  matrix form
\begin{align}\label{eq:ADMM_eq_matrix_x}
\begin{bmatrix}
x^{k+1}\\
y^{k}
\end{bmatrix}
&=
%
\underbrace{\begin{bmatrix}
M_{11}    &  M_{12}    \\
M_{21} & (1-\alpha)I
\end{bmatrix}
}_M
\begin{bmatrix}
x^{k}\\
y^{k-1}
\end{bmatrix}
,
\end{align}
for $k\geq 1$ with $x^1 = -(Q+\rho\bar{E}^\top\bar{E})^{-1}\left(q+\rho\bar{E}^\top (\bar{F}z^0+u^0)\right)$, $y^0=\bar{E}^\top \bar{F}z^0$, $z^0=-(\bar{F}^\top\bar{F})^{-1}c/\rho$, $u^0=\bar{F}z^0$, and
\begin{equation}
\label{eq:ADMM_M11_M12}
\begin{aligned}
M_{11} &=\alpha\rho(Q+\rho \bar{E}^\top \bar{E})^{-1}\bar{E}^\top\left(\Pi_{\mathcal{R}(\bar{F})} -  \Pi_{\mathcal{N}(\bar{F}^\top)}\right)\bar{E} + I,\\
M_{12}&=   \alpha\rho(Q+\rho \bar{E}^\top \bar{E})^{-1},\quad
M_{21}=  -\alpha \bar{E}^\top \Pi_{\mathcal{R}(\bar{F})} \bar{E}.
\end{aligned}
\end{equation}

The next theorem shows how the convergence properties of the ADMM iterations are characterized by the spectral properties of the matrix $M$.
\begin{thm}\label{thm:convergence_ADMM}
Define $\sigma^{k+1} \triangleq [x^{{k+1}^\top}\; y^{{k}^\top}]^\top$, $s\triangleq\mbox{dim}\left(\mathcal{R}(\bar{F})\cap \mathcal{R}(\bar{E})\right)$, and let $\{\phi_i\}$ be the eigenvalues of $M$ ordered so that $|\phi_1| \leq \dots \leq \dots \leq |\phi_{2n}|$. The ADMM iterations~\eqref{eqn:admm_scaled_iterations} converge to the optimal solution of~\eqref{eq:QP_problem} if and only if $s\geq1$ and $1=\phi_{2n} = \dots = \phi_{2n-s+1} > |\phi_{2n-s}|$.  Moreover, the convergence factor of the ADMM iterates in terms of the sequence $\{\sigma^k\}$ equals  $\phi^\star = |\phi_{2n-s}|$.
\end{thm}
\begin{IEEEproof}
See Appendix~\ref{app:thm:convergence_ADMM}.
\end{IEEEproof}
Below we state the main problem to be addressed in the remainder of this paper.
\begin{problem}\label{prob:optimal_scaling}
Which scalars $\rho^\star$ and $\alpha^\star$ and what matrix $R^\star$ minimize $\vert\phi_{2n-s}\vert$, the convergence factor of the ADMM iterates?
\end{problem}

As the initial step to tackle Problem~\ref{prob:optimal_scaling}, we characterize the eigenvalues $\phi_i$ of $M$.
Our analysis will be simplified by choosing an $R$ that satisfies the following assumption.
\begin{assumption}\label{assum:EWE_Q}
The scaling matrix $R$ is such that $E^\top R^\top R E = \bar{E}^\top \bar{E} = \kappa Q$ for some $\kappa>0$ and $\bar{E}^\top \bar{E} \succ 0$.
\end{assumption}

Assumption~\ref{assum:EWE_Q} may appear restrictive at first sight, but we will later describe several techniques for finding such an $R$, even for the distributed setting outlined in Section~\ref{sec:background}. 
Replacing $\bar{E}^\top \bar{E} = \kappa Q$ in~\eqref{eq:ADMM_M11_M12} and using the identity $\Pi_{\mathcal{R}(\bar{F})} -  \Pi_{\mathcal{N}(\bar{F}^\top)}= 2\Pi_{\mathcal{R}(\bar{F})} -  I$ yields
\begin{equation*}
\begin{aligned}
M_{11} &=\alpha\dfrac{\rho\kappa}{1+\rho \kappa}(\bar{E}^\top \bar{E})^{-1}\bar{E}^\top\left(2\Pi_{\mathcal{R}(\bar{F})} -  I\right)\bar{E} + I,\\
M_{12}&=  \alpha\dfrac{\rho\kappa}{1+\rho \kappa}(\bar{E}^\top \bar{E})^{-1},\quad
M_{21}=  -\alpha \bar{E}^\top \Pi_{\mathcal{R}(\bar{F})} \bar{E}.
\end{aligned}
\end{equation*}
These expressions allow us to explicitly characterize the eigenvalues of $M$~in~\eqref{eq:ADMM_eq_matrix_x}.

\begin{thm}\label{thm:M_eigenvalues}
Consider the ADMM iterations~\eqref{eq:ADMM_eq_matrix_x} and suppose that $\bar{E}^\top \bar{E} = \kappa Q$. Let $v_i$ be a generalized eigenvector of $\left(\bar{E}^\top\left(2\Pi_{\mathcal{R}(\bar{F})} - I\right)\bar{E},\;  \bar{E}^\top \bar{E}\right)$ with associated generalized eigenvalue $\lambda_i$. Then, $M$ has two right eigenvectors on the form $[v_i^\top \; w_{i1}^\top]^\top$ and $[v_i^\top \; w_{i2}^\top]^\top$ whose associated eigenvalues $\phi_{i1}$ and $\phi_{i2}$ are the solutions to the quadratic equation

\begin{equation}\label{eq:eigenvalues_M_phi}
  \begin{aligned}
    \phi_i^2+ a_1(\lambda_i)\phi_i + a_0(\lambda_i)=0,
  \end{aligned}
\end{equation}
where
\begin{equation}\label{eq:polynomial_coefficients}
  \begin{aligned}
     a_1(\lambda_i)&\triangleq  \alpha - \alpha \beta{\lambda_i} - 2 , \quad \beta \triangleq {\rho \kappa} / {(1+\rho \kappa)},\\
     a_0(\lambda_i)&\triangleq \alpha \beta(1-\dfrac{\alpha}{2}) \lambda_i + \dfrac{1}{2}\alpha^2 \beta+1-\alpha.
  \end{aligned}
\end{equation}
\begin{IEEEproof}
See Appendix~\ref{app:thm:M_eigenvalues}.
\end{IEEEproof}
\end{thm}

From~\eqref{eq:eigenvalues_M_phi} and \eqref{eq:polynomial_coefficients} one directly sees that $\alpha$, $\rho$ (or, equivalently, $\beta$) and $R$ affect the eigenvalues of $M$. We will use $\phi(\alpha, \beta, \lambda_i)$ to emphasize this dependence.
%
In the next section we study the properties of~\eqref{eq:eigenvalues_M_phi} with respect to $\beta$, $\alpha$, and $\lambda_i$. 

%
\subsection{Optimal parameter selection}
To minimize the convergence factor of the iterates \eqref{eqn:admm_scaled_iterations}, we combine Theorem~\ref{thm:convergence_ADMM}, which relates the convergence factor of the ADMM iterates to the spectral properties of the matrix $M$, with Theorem~\ref{thm:M_eigenvalues}, which gives explicit expressions for the eigenvalues of $M$ in terms of the ADMM parameters. The following result useful for the development of our analysis.
\begin{proposition}[Jury's stability test~\cite{Jury_stability_test}]\label{prop:Jury}
The quadratic polynomial $a_2 \phi_i^2 + a_1 \phi_i + a_0$ with real coefficients $a_2> 0$, $a_1$, and $a_0$ has its roots inside the unit-circle, \ie, $|\phi_i|<1$, if and only if the following three conditions hold:
\begin{equation*}
\begin{aligned}
\mbox{i)}\quad   &  a_0 + a_1 + a_2 > 0;\\
\mbox{ii)}\quad  & a_2 > a_0;\\
\mbox{iii)}\quad & a_0 - a_1 + a_2 > 0.
\end{aligned}
\end{equation*}
\end{proposition}

The next sequence of lemmas derive some useful properties of $\lambda_i$ and of the eigenvalues of $M$.

\begin{lem}\label{lem:bar_lambda}
The generalized eigenvalues of $(\bar{E}^\top\left(2\Pi_{\mathcal{R}(\bar{F})} - I\right)\bar{E},\, \bar{E}^\top \bar{E})$ are real scalars in $[-1,1]$.
\end{lem}
\begin{IEEEproof}
See Appendix~\ref{app:lem:bar_lambda}.
\end{IEEEproof}

\begin{lem}\label{lem:bar_lambda_ones}
Let $\lambda_i$ be the $i$-th generalized eigenvalue of $(\bar{E}^\top\left(2\Pi_{\mathcal{R}(\bar{F})} - I\right)\bar{E},\, \bar{E}^\top \bar{E})$, ordered as $\lambda_n \geq \dots \geq\lambda_i \geq \dots\geq \lambda_1$ and let $\dim \left(\mathcal{R}(\bar{E}) \cap \mathcal{R}(\bar{F}) \right)= s$. If the optimization problem~\eqref{eq:QP_problem} is feasible, we have $s \geq 1$ and $\lambda_i=1$, for all $i = n, \dots, n-s+1$.
\end{lem}
\begin{IEEEproof}
See Appendix~\ref{app:lem:bar_lambda_ones}.
\end{IEEEproof}

\begin{lem}
\label{lem:admm_stability}
Consider the eigenvalues $\{\phi_i\}$ of the matrix $M$ in~\eqref{eq:ADMM_eq_matrix_x}, ordered as $|\phi_{2n}| \geq \dots \geq |\phi_{i}|\geq \dots \geq |\phi_{1}|$. It follows that $\phi_{2n} = \dots= \phi_{2n-s+1} = 1$ where $s=\dim \left(\mathcal{R}(\bar{E}) \cap \mathcal{R}(\bar{F}) \right)$. Moreover, for $\beta \in (0,\, 1)$ and $\alpha\in (0,\,2]$ we have $|\phi_i|< 1 $ for $i\leq 2n-s$.
\end{lem}
\begin{IEEEproof}
See Appendix~\ref{app:lem:admm_stability}.
\end{IEEEproof}

Lemma~\ref{lem:admm_stability} and Theorem~\ref{thm:convergence_ADMM} establish that the convergence factor of the ADMM iterates, $  |\phi_{2n-s}|$, is strictly less than $1$ for $\beta \in(0,\,1)$ and $\alpha\in(0,\,2]$.
 Next, we characterize $\vert \phi_{2n-s}\vert$ explicitly in terms of $\alpha$, $\beta$ and $\lambda_i$.
\begin{thm}\label{thm:reduced_max}
Consider the eigenvalues $\{\phi_i\}$ of $M$ ordered as in Lemma \ref{lem:admm_stability}.
For fixed $\alpha\in(0,\,2]$ and $\beta\in(0,\,1)$, the magnitude of $\phi_{2n-s}$ is given by
\begin{align}
\label{eq:convergence_factor_M}
%
\vert \phi_{2n-s}\vert \triangleq \max\left\{ g^{+}_r,\; g^{-}_r,\; g_c,\; g_1 \right\}
\end{align}
where
\begin{equation}\label{eq:function_g_lambda}
\begin{aligned}
 g^+_r &\triangleq 1 + \frac{\alpha}{2}\beta\lambda_{n-s} - \frac{\alpha}{2} + \frac{\alpha}{2}\sqrt{\lambda_{n-s}^2 \beta^2 - 2\beta +1 + s^+_r},\\
 g^-_r &\triangleq  - 1 - \frac{\alpha}{2}\beta\lambda_1 + \frac{\alpha}{2} + \frac{\alpha}{2}\sqrt{\lambda_1^2 \beta^2 - 2\beta +1 + s^{-}_r},\\
 g_c &\triangleq \sqrt{\frac{1}{2} \alpha^2\beta(1-\lambda_{n-s}) +1 -\alpha + \alpha \beta\lambda_{n-s} + s_c},\\
 g_1 & \triangleq |1-\alpha(1-\beta)|,\\
s^+_r &\triangleq \max\{0,\; - (\beta^2\lambda_{n-s}^2 -2\beta +1  ) \},  \\
s^{-}_r &\triangleq  \max\{0,\; - (\beta^2\lambda_{1}^2 -2\beta +1  ) \}, \\
s_c &\triangleq  \max\{0,\; - a_0(\lambda_{n-s})  \}.
\end{aligned}
\end{equation}

Moreover, we have $\vert \phi_{2n-s}\vert > g^+_r$, $\vert \phi_{2n-s}\vert > g^-_r$, and $\vert \phi_{2n-s}\vert > g_c$ if $s^+_r > 0$, $s^-_r > 0$, and $s_c  > 0$, respectively.
\end{thm}
\begin{IEEEproof}
see Appendix~\ref{app:thm:reduced_max}.
\end{IEEEproof}

Given the latter result, the problem of minimizing $|\phi_{2n-s}|$ with respect to $\alpha$ and $\beta$ can be written as
\begin{equation*}
\begin{aligned}
\min_{\alpha\in(0,\, 2],\; \beta \in(0,\,1)} \max\left\{ g^+_r,\; g^-_r,\; g_c ,\; g_1  \right\}.
\end{aligned}
\end{equation*}
Numerical studies have suggested that under-relaxation, \ie, letting $\alpha < 1$, does not improve the convergence speed of ADMM, see, \eg, \cite{Boyd11}.
The next result establishes formally that this is indeed the case for our considered class of problems.
\begin{proposition}\label{prop:alpha_larger_than_one}
Let $\beta\in(0,\,1)$ be fixed and consider $\phi_{2n-s}(\alpha,\beta)$ . For $\alpha<1$, it holds that $\vert \phi_{2n-s}(1,\beta) \vert < \vert \phi_{2n-s}(\alpha,\beta) \vert$.
\end{proposition}
\begin{IEEEproof}
See Appendix~\ref{app:prop:alpha_larger_than_one}.
\end{IEEEproof}

The main result presented below provides explicit expressions for the optimal parameters $\alpha$ and $\beta$ that minimize $\vert \phi_{2n-s}\vert$ over given intervals.
\begin{thm}\label{thm:optimal_f_alpha}
Consider the optimization problem~\eqref{eq:QP_problem} under Assumption~\ref{assum:EWE_Q} and its associated ADMM iterates~(\ref{eq:ADMM_eq_matrix_x}). The parameters $\alpha^\star$ and $\beta^\star$ that minimize the convergence factor $|\phi_{2n-s}|$ over $\alpha\in(0,\alpha^\star]$ and $\beta\in(0,1)$ are:
\begin{enumerate}[leftmargin=*,widest=III, align=right, label=Case \Roman*]
\item: if $\lambda_{n-s}>0$ and $\lambda_{n-s}\geq |\lambda_1|$,
\begin{equation}
\begin{aligned}
\beta^\star &= \dfrac{1-\sqrt{1-\lambda_{n-s}^2}}{\lambda_{n-s}^2}, \quad
\alpha^\star = 2,\\
|\phi_{2n-s}| &= \dfrac{1-\sqrt{1-\lambda_{n-s}^2}}{\lambda_{n-s}};
\end{aligned}
\end{equation}\label{enum:case1}
\item: if $|\lambda_1|\geq\lambda_{n-s}>0$,
\begin{equation}\label{eq:alpha_beta_suboptimal_1}
\begin{aligned}
\beta^\star &= \dfrac{1-\sqrt{1-\lambda_{n-s}^2}}{\lambda_{n-s}^2},\\
\alpha^\star &=
\dfrac{4}{2-\left(\lambda_{n-s} + \lambda_1 - \sqrt{\lambda_1^2 - \lambda_{n-s}^2}\right)\beta^\star},\\
|\phi_{2n-s}| &= 1 + \frac{\alpha^\star}{2}\lambda_{n-s}\beta^\star - \frac{\alpha^\star}{2};
\end{aligned}
\end{equation}\label{enum:case2}
\item: if $0\geq\lambda_{n-s}\geq\lambda_1$,
\begin{equation}
\begin{aligned}
\beta^\star &= \dfrac{1}{2}, \quad
\alpha^\star = \dfrac{4}{2- \lambda_1},\quad
|\phi_{2n-s}| =  \dfrac{-\lambda_1}{2-\lambda_1}.
\end{aligned}
\end{equation}\label{enum:case3}
\end{enumerate}
\end{thm}
\begin{IEEEproof}
  See Appendix~\ref{app:thm:optimal_f_alpha}.
\end{IEEEproof}


Considering the standard ADMM iterations with $\alpha=1$, the next result immediately follows.
\begin{cor}\label{cor:optimal_parameters_alpha_1}
For $\alpha=1$, the optimal $\beta^\star$ that minimizes the convergence factor $\vert \phi_{2n-s}^\star\vert$ is
 \begin{align}\label{eq:optimal_beta_alpha_1}
  \beta^\star =\left\{
  \begin{array}[c]{ll}
   \dfrac{1-\sqrt{1-\lambda^2_{n-s}}}{\lambda_{n-s}^2}& \lambda_{n-s}> 0,\\
   \dfrac{1}{2}& \lambda_{n-s}\leq 0.
  \end{array}\right.
  \end{align}
  Moreover, the corresponding convergence factor is
  \begin{align}\label{eq:optimal_factor_alpha_1}
  \vert \phi_{2n-s}^\star\vert = \left\{
  \begin{array}[c]{ll}
   \dfrac{1}{2}\left(1+\dfrac{\lambda_{n-s}}{1+\sqrt{1-\lambda_{n-s}^2}}\right)& \lambda_{n-s}> 0,\\
   \dfrac{1}{2}& \lambda_{n-s}\leq 0.
  \end{array}\right.
  \end{align}
\end{cor}
\begin{IEEEproof}
From the proof of Theorem~\ref{thm:optimal_f_alpha}, when $\lambda_{n-s}>0$, then $\beta^\star = {(1-\sqrt{1-\lambda_{n-s}^2})}/{\lambda_{n-s}^2}$ is optimal, and when $\lambda_{n-s}\leq 0$, $\beta^\star = 1/2$ is the minimizer. The result follows by setting $\alpha=1$ and obtaining corresponding convergence factors that are given by $g_r^{+}(\alpha=1, \beta^\star, \lambda_{n-s})$.
\end{IEEEproof}

\subsection{Optimal constraint scaling}\label{sec:Optimal_const_scaling}
As seen in Theorem~\ref{thm:optimal_f_alpha}, the convergence factor of ADMM depends in a piecewise fashion on $\lambda_{n-s}$ and $\lambda_{1}$. In the first two cases, the convergence factor is monotonically increasing in $\lambda_{n-s}$, and it makes sense to choose the constraint scaling matrix $R$ to minimize $\lambda_{n-s}$ while satisfying the structural constraint imposed by Assumption~\ref{assum:EWE_Q}.
To formulate the selection of $R$ as a quasi-convex optimization problem, we first enforce the constraint $\kappa Q = E^\top W E$ by using the following result.
\begin{lem}\label{lem:EWE_Qkappa_general}
Consider the optimization problem~\eqref{eq:QP_problem} with $Q\succeq0$ and let $P\in\mathbf{R}^{n\times s}$ be an orthonormal basis for $\Null{\Pi_{\mathcal{N}(F^\top)}E}$. Let $W=R^{\top}R$ and assume that $E^\top W E \succ 0$. If $P^\top E^\top W E P = P^\top Q P\succ 0$, then the optimal solution to~\eqref{eq:QP_problem} remains unchanged when $Q$ is replaced with $E^\top W E$.
\end{lem}
\begin{IEEEproof}
See Appendix~\ref{app:lem:EWE_Qkappa_general}.
\end{IEEEproof}

The following result addresses Assumption~\ref{assum:feasible}.
\begin{lem}\label{lem:feasible}
Let $P_1$ be an orthonormal basis for the orthogonal complement to $\Null{\Pi_{\mathcal{N}(F^\top)}E}$ and define $W=R^\top R\succeq 0$. The following statements are true:
 \begin{enumerate}[label=\roman*)]
 \item Assumption~\ref{assum:feasible} holds if and only if $F^\top W F\succ0$ and $P_1^\top\bar{E}^\top\Pi_{\Null{\bar{F}^\top}}\bar{E}P_1 \succ0$;
 \item If Assumption~\ref{assum:feasible} holds, then $\Null{\Pi_{\mathcal{N}(F^\top)}E} = \Null{\Pi_{\mathcal{N}(\bar{F}^\top)}\bar{E}}$.
 \end{enumerate}
\end{lem}
\begin{IEEEproof}
See Appendix~\ref{app:lem:feasible}.
\end{IEEEproof}

Next, we derive a tight upper bound on $\lambda_{n-s}$.
\begin{lem}\label{lem:LMI}
Let $P_1$ be an orthonormal basis for the orthogonal complement to $\Null{\Pi_{\mathcal{N}(F^\top)}E}$. Defining $W=R^\top R\succeq 0$ and letting $\lambda\leq 1$, we have $\lambda \geq \lambda_{n-s}$ if and only if
\begin{equation}\label{eq:upperbound_lambda}
 \begin{aligned}
\begin{bmatrix}
(\lambda+1) P_1^\top E^\top W E  P_1 &  P_1^\top E^\top WF \\
F^\top W E  P_1               &   \dfrac{1}{2} F^\top W F
\end{bmatrix}
 \succ 0.
\end{aligned}
\end{equation}
Moreover, Assumption~\ref{assum:feasible} holds for a given $W$ satisfying~\eqref{eq:upperbound_lambda} with $\lambda\leq 1$.
\end{lem}
\begin{IEEEproof}
See Appendix~\ref{app:lem:LMI}.
\end{IEEEproof}

Using the previous results, the matrix $W$ minimizing $\lambda_{n-s}$ can be computed as follows.
\begin{thm}\label{thm:weight_optimization_general}
Let $P_1\in\mathbf{R}^{n\times n-s}$ be an orthonormal basis for the orthogonal complement to $\Null{\Pi_{\mathcal{N}(F^\top)}E}$, define $P\in\mathbf{R}^{n\times s}$ as an orthonormal basis for $\Null{\Pi_{\mathcal{N}(F^\top)}E}$,  and denote $\mathcal{A}$ as a given sparsity pattern. The matrix $W=R^\top R \in \mathcal{A}$ that minimizes
$\lambda_{n-s}$ while satisfying Assumptions~\ref{assum:EWE_Q} and~\ref{assum:feasible} is the solution to the quasi-convex optimization problem
\begin{equation}\label{eqn:weight_optimization_general}
\begin{array}{ll}
\underset{W,\lambda}{\mbox{minimize}} & \lambda\\
\mbox{subject to} & W\in\mathcal{A},\; W \succeq 0,\;\eqref{eq:upperbound_lambda},\\
 & P^\top E^\top W E P = P^\top Q P.

\end{array}
\end{equation}
\end{thm}
\begin{IEEEproof}
The proof follows from Lemmas~\ref{lem:EWE_Qkappa_general},~\ref{lem:feasible}, and~\ref{lem:LMI}. 
\end{IEEEproof}

The results derived in the present section contribute to improve the convergence properties of the ADMM algorithm for equality-constrained quadratic programming problems. The procedure to determine suitable choices of $\rho$, $\alpha$, and $R$ is summarized in Algorithm~\ref{alg:optimal_scaling_general}.
\begin{algorithm}[]                      
\caption{Optimal Constraint Scaling and Parameter Selection}          
\begin{enumerate}                    
\item Compute $W^\star$ and the corresponding $\lambda_{n-s}$ and $\lambda_{1}$ according to Theorem~\ref{thm:weight_optimization_general};
\item Using Lemma~\ref{lem:EWE_Qkappa_general}, replace $Q$ with $E^\top W E$ and let $\kappa=1$;
\item Given $\lambda_{n-s}$ and $\lambda_1$, use the ADMM parameters $\rho^\star=\frac{\beta^\star}{1-\beta^\star}$ and $\alpha^\star$ proposed in Theorem~\ref{thm:optimal_f_alpha}.
\end{enumerate}
\label{alg:optimal_scaling_general}
\end{algorithm}

\section{ADMM for distributed quadratic programming}\label{sec:distributed_QP}
We are now ready to develop optimal scalings for the ADMM iterations for distributed quadratic programming. Specifically, we consider~\eqref{eq:DQP_consensus}
with $f_i(x_{i}) = (1/2) x_i^\top Q_i x_i + q_i^\top x_i$ and $Q_i\succ 0$ and use the results derived in the previous section to derive optimal algorithm parameters for the ADMM iterations in both edge- and node-variable formulations.

\subsection{Enforcing agreement with edge variables}\label{subsec:Edge_variable_detail}
In the edge variable formulation, we introduce  auxiliary variables $z_{\{i,j\}}$ for each edge $\{i,j\}\in\mathcal{E}$ and re-write the optimization problem in the form of (\ref{eq:DQP_high}).
The resulting ADMM iterations for node $i$ can be written as
\begin{equation}\label{eq:ADMM_edge_i}
\begin{aligned}
x_i^{k+1} &= \underset{x_i}{\mbox{argmin}}\; \frac{1}{2} x_i^\top Q_i x_i + q_i^\top x_i
\iftoggle{draft}{%
%
}{%
  \\&\quad\quad
}
+  \frac{\rho}{2} \sum_{j\in\mathcal{N}_i} \|R_{(i,j)} x_i - R_{(i,j)} z_{\{i,j\}}^k  + R_{(i,j)} u_{(i,j)}^{k} \|_2^2, \\
\gamma_{(j,i)}^{k+1} &= \alpha x_j^{k+1} + (1-\alpha)z_{\{i,j\}}^k,\quad \forall j\in\mathcal{N}_i ,\\
%
z_{\{i,j\}}^{k+1} &= \underset{z_{\{i,j\}}}{\mbox{argmin}}\;   \|R_{(i,j)} \gamma_{(i,j)}^{k+1} + R_{(i,j)} u_{(i,j)}^{k}   - R_{(i,j)} z_{\{i,j\}}  \|_2^2
\iftoggle{draft}{%
%
}{%
  \\
&\quad\qquad\,
}
+  \|R_{(j,i)}\gamma_{(j,i)}^{k+1}+ R_{(j,i)} u_{(j,i)}^{k} - R_{(j,i)} z_{\{i,j\}}  \|_2^2,\\
u_{(i,j)}^{k+1} &=u_{(i,j)}^{k} + \gamma_{(i,j)}^{k+1} - z_{\{i,j\}}^{k+1}. \\
\end{aligned}
\end{equation}
Here, $u_{(i,j)}$ is the scaled Lagrange multiplier, private to node $i$, associated with the constraint $R_{(i,j)} x_i = R_{(i,j)} z_{\{i,j\}}$, and the variables $\gamma_{(i,j)}$ have been introduced to write the iterations in a more compact form. Note that the algorithm is indeed distributed, since each node $i$ only needs the current iterates $x_j^{k+1}$ and $u_{(j,i)}^{k}$ from its neighboring nodes $j\in\mathcal{N}_i$.

We can also re-write the problem formulation as an equality constrained quadratic program on the form (\ref{eq:consensus_problem_edge}) with $f(x) = (1/2)x^\top Q x + q^\top x$ , $Q=\mbox{diag}\left(\left\{{Q}_i\right\}_{i\in\mathcal{V}}\right)$, and $q^\top = [q_1^\top\,\dots\, q_{\vert\mathcal{V}\vert}^\top]$.
As shown in Section~\ref{sec:qp_equality}, the associated ADMM iterations can be written in vector form (\ref{eqn:admm_scaled_iterations}) and the step-size and the relaxation parameter that minimize the convergence factor of the iterates are given in Theorem~\ref{thm:optimal_f_alpha}.

Recall the assumptions that $W\succeq 0$ is chosen so that $E^\top W E = \kappa Q$ for $\kappa>0$. The next result shows that such assumptions can be satisfied locally by each node.
\begin{lem}\label{lem:local_weights}
Consider the distributed optimization problem described by~\eqref{eq:consensus_problem_edge} and~\eqref{eq:consensus_problem_edge_matrices} and let $W=R^\top R$. The equation $E^\top W E=\kappa Q$ can be ensured for any $\kappa>0$ by following a weight-assignment scheme satisfying the local constraints $\sum_{j\in \mathcal{N}_i} W_{(i,j)} = \kappa Q_i$ for all $i\in\mathcal{V}$.
\end{lem}
\begin{IEEEproof}
From the $x_i$-update in the ADMM iterations~\eqref{eq:ADMM_edge_i}, we see that the diagonal block of $E^\top W E$ corresponding to node $i$ is given by $\sum_{j\in \mathcal{N}_i} W_{(i,j)}$. Hence, $E^\top W E=\kappa Q$ is met if each agent $i$ ensures that $\sum_{j\in \mathcal{N}_i} W_{(i,j)} = \kappa Q_i$.
\end{IEEEproof}

Next, we analyze in more detail the scalar case with symmetric edge weights.
\subsubsection{Scalar case}
Consider the scalar case $n_x =1$ with $n=\vert \mathcal{V} \vert$ and let the edge weights be symmetric with $W_{(i,j)}=W_{(j,i)}=w_{\{i,j\}} \geq 0$ for all $(i,j)\in\bar{\mathcal{E}}$.
As derived in Section~\ref{sec:qp_equality}, the ADMM iterations can be written in matrix form as~\eqref{eq:ADMM_eq_matrix_x}. Exploiting the structure of $E$ and $F$, we derive
\begin{equation*}
\begin{aligned}
M_{11} &= \alpha\rho(Q+\rho D)^{-1}A + I,\quad   M_{12} =   \alpha \rho(Q+\rho D)^{-1},
\iftoggle{draft}{%
\quad  M_{21} = -\dfrac{\alpha}{2}(D+A).
}{%
\\ M_{21} &= -\dfrac{\alpha}{2}(D+A).
}
\end{aligned}
\end{equation*}

The optimal step-size $\rho^\star$ and $\alpha^\star$ that minimizes the convergence factor $\vert\phi_{2n-1}\vert$ are given in Theorem~\ref{thm:optimal_f_alpha}, where the eigenvalues $\{\lambda_i\}$ in the corresponding theorems are the set of ordered generalized eigenvalues of $(A,D)$. Here we briefly comment on the relationship between the generalized eigenvalues of $(A,D)$ and the eigenvalues of the normalized Laplacian, $L_{\text{norm}}= I - D^{-1/2}A D^{-1/2}$. In particular, we have $1 - \lambda= \psi$,
where $\psi$ is any eigenvalue of the normalized Laplacian and $\lambda$ is a generalized eigenvalue of $(A,D)$ corresponding to $\psi$. For certain well-known classes of graphs the value of the eigenvalues of the normalised Laplacian are known, e.g. see \cite{chung1997spectral} and \cite{butler2008eigenvalues} for more information. As a result, one can identify  which case of Theorem~\ref{thm:optimal_f_alpha} is  applied to each of these graphs. The following proposition establishes one such result.
\begin{proposition}
Adopt the hypothesis of Theorem~\ref{thm:optimal_f_alpha}. The following statements are true.
\begin{enumerate}[label=\roman*)]
\item Case II of Theorem~\ref{thm:optimal_f_alpha} holds for path graphs with $|\mathcal{V}|\geq 4$, cycle graphs with $|\mathcal{V}| \geq 5$, and wheel-graphs with $|\mathcal{V}|\geq 6$.
\item Case III of Theorem~\ref{thm:optimal_f_alpha} holds for complete graphs,  bi-partite graphs,  star graphs, path graphs with $|\mathcal{V}|=3$, cycle graphs with $|\mathcal{V}|\in\{3,4\}$ and wheel-graphs with $|\mathcal{V}|\in\{4,5\}$.
\end{enumerate}
\end{proposition}
\begin{IEEEproof}
The proof is a direct consequence of the analytical expressions of the eigenvalues of the normalised Laplacian of given in \cite{chung1997spectral,butler2008eigenvalues} and the relationship $1 - \lambda= \psi$.
\end{IEEEproof}

For general topologies, without computing the generalized eigenvalues it is not easy to know which case of Theorem~\ref{thm:optimal_f_alpha} applies. Moreover, when we use non-unity edge-weights, optimizing these for one case might alter the generalized eigenvalues so that another case applies. In extensive simulations, we have found that a good heuristic is to use scalings that attempt to reduce the magnitude of both the smallest and the second-largest generalized eigenvalues. The next lemma shows how to compute such scalings for the edge-variable formulation.
\begin{lem}\label{lem:weight_optimization_edge}
Consider the weighted undirected graph $\mathcal{G}=(\mathcal{V},\mathcal{E},\mathcal{W})$. The non-negative edge-weights $\{w_{\{i,j\}}\}$ that jointly minimize and maximize the second largest and smallest generalized eigenvalue of $(A,D)$, $\lambda_{n-1}$ and $\lambda_1$, are obtained from the optimal solution to the quasi-convex problem
\begin{equation}
\label{eq:weight_optimization_edge}
\begin{array}[c]{lll}
\underset{\{w_{\{i,j\}}\},\,\lambda}{\mbox{minimize}} & \lambda &\\
\mbox{subject to}
%
& w_{\{i,j\}} \geq 0,   & \forall \,i,j\in\mathcal{V},\\
& A_{ij} = w_{\{i,j\}},  & \forall \,\{i,j\}\in\mathcal{E},\\
& A_{ij} = 0,                &\forall \,\{i,j\}\not\in\mathcal{E},\\
& D = \mbox{diag}(A\textbf{1}_n),&\\
& D \succ \epsilon I,&\\
%
& P^\top\left( A - \lambda D \right) P \prec 0,&\\
& A + \lambda D \succ 0, &
\end{array}
\end{equation}
where the columns of $P\in\mathbf{R}^{n\times n-1}$ form an orthonormal basis of $\mathcal{N}(\textbf{1}_n^\top)$ and $\epsilon>0$.
\begin{IEEEproof}
The second last constraint ensures that $\lambda > \lambda_{n-1}$ and follows from a special case of Lemma~\ref{lem:LMI}, while the last constraint enforces $\lambda_1 > - \lambda$.
\end{IEEEproof}
\end{lem}

\subsection{Enforcing agreement with node variables}\label{subsec:Node_variable_detail}
Recall the node variable formulation~\eqref{eq:DQP_node} using the auxiliary variables $z_i\in\mathbf{R}^{n_x}$ for each node $i\in\mathcal{V}$ described in Section~\ref{subsec:Node_variable}.
The ADMM iterations for node $i$ can be rewritten as
\begin{equation}\label{eq:ADMM_node_i}
\begin{aligned}
x_i^{k+1} &= \underset{x_i}{\mbox{argmin}}\: \frac{1}{2} x_i^\top Q_i x_i + q_i^\top x_i
\iftoggle{draft}{%
%
}{%
  \\
&\quad\quad
}
+\frac{\rho}{2}\sum_{j\in\mathcal{N}_i\cup\{i\}} \!\! \|R_{(i,j)} x_i - R_{(i,j)} z_{j}^k  + R_{(i,j)} u_{(i,j)}^{k} \|_2^2, \\
\gamma_{(j,i)}^{k+1} &= \alpha x_j^{k+1} + (1-\alpha)z_{i}^k, \qquad\forall j\in\mathcal{N}_i\cup\{i\},\\
z_{i}^{k+1} &= \underset{z_{i}}{\mbox{argmin}} \! \!\!
\sum_{j\in \mathcal{N}_i\cup\{i\}} \!\!\! \|R_{(j,i)} \gamma_{(j,i)}^{k+1} + R_{(j,i)} u_{(j,i)}^{k}   - R_{(j,i)} z_{i}  \|_2^2,\\
%
u_{(i,j)}^{k+1} &=u_{(i,j)}^{k} + \gamma_{(i,j)}^{k+1} - z_{j}^{k+1},  \qquad\forall j\in\mathcal{N}_i\cup\{i\},\\
\end{aligned}
\end{equation}
where $u_{(i,j)}$ is the scaled Lagrange multiplier, private to node $i$, associated with the constraint $R_{(i,j)} x_i = R_{(i,j)} z_{j}$, and $\gamma_{(i,j)}(k)$ is an auxiliary variable private to node $i$ and associated with the edge $(i,j)$. Note that the algorithm is distributed, since it only requires communication between neighbors. However, unlike the previous formulation with edge variables, here two communication rounds must take place: the first to exchange the private variables $x_j^{k+1}$ and $u_{(j,i)}^{k}$, required for the $z_{i}$-update; the second to exchange the private variables $z_j^{k+1}$, required for the $u_{(i,j)}$- and $x_i$-updates.

Let $x=[x_1^\top\,\cdots\, x_{\vert \mathcal{V}\vert}^\top]^\top$, $z = [z_{1}^\top\, \cdots \, z_{\vert \mathcal{V}\vert}^\top]^\top$, $Q=\mbox{diag}\left( \left\{{Q}_i \right\}_{i\in\mathcal{V}}\right)$, and $q^\top = [{q}_1\,\dots\, {q}_{\vert\mathcal{V}\vert}]$. The cost function in~\eqref{eq:DQP_node} takes the form $f(x) = (1/2)x^\top Q x + q^\top x$ while $E$ and $F$ are given in~\eqref{eq:consensus_problem_node_matrices}.
Note that the matrix $E$ is the same as for the edge-variable case, thus Lemma~\ref{lem:local_weights} may also be applied to the present formulation to ensure $ E^\top W E = \kappa Q$.


Next we consider the scalar case with symmetric weights.

\subsubsection{Scalar case}
Consider the scalar case $n_x =1$ with $n=\vert \mathcal{V} \vert$ and let the edges weights be symmetric with $W_{(i,j)}=W_{(j,i)}=w_{\{i,j\}}\geq0$ for all $(i,j)\in\bar{\mathcal{E}}$.
Using the structure of $E$ and $F$, the fixed point equation~\eqref{eq:ADMM_eq_matrix_x} can be formulated by the following relations
\begin{equation}\label{eq:Consensus_node_linearMatrix}
\begin{aligned}
 M_{11} &= \alpha\rho(Q+\rho D)^{-1}(2AD^{-1}A-D)+ I,
 \iftoggle{draft}{%
\quad   M_{12} =   \alpha\rho(Q+\rho D)^{-1}, \quad  M_{21} = -\alpha AD^{-1}A.
}{%
\\  M_{12} &=   \alpha\rho(Q+\rho D)^{-1}, \quad M_{21} = -\alpha AD^{-1}A.
}
\end{aligned}
\end{equation}

The optimal step-size $\rho^\star$ and $\alpha^\star$ that minimizes the convergence factor $\vert\phi_{2n-s}\vert$ are given in Theorem~\ref{thm:optimal_f_alpha}, where the eigenvalues $\{\lambda_i\}$ are the set of ordered generalized eigenvalues of $(2AD^{-1}A-D,\;D)$.
For the node-variable formulation, we have not been able to formulate a weight optimization corresponding to Lemma~\ref{lem:weight_optimization_edge}. However, in the numerical evaluations, we will propose a modification that ensures that Case I of Theorem~\ref{thm:optimal_f_alpha} applies, and then minimize the second-largest generalized eigenvalue.
\section{Numerical examples}
\label{sec:numerical}
Next, we illustrate our results via numerical examples.
\subsection{Distributed quadratic programming}
\begin{figure}[htbp]
  \centering
  \includegraphics[width=0.6\hsize]{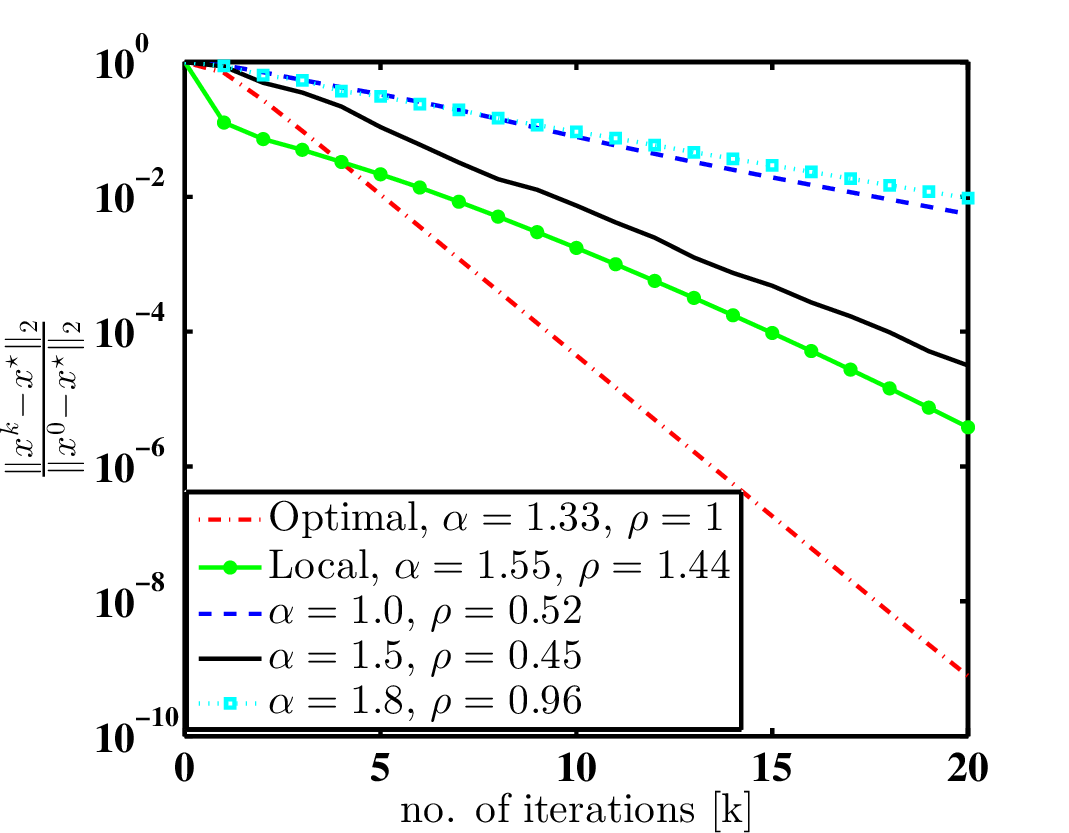}
  \vspace{-10pt}
\caption{Normalized error for the scaled ADMM algorithm with $W^\star$ from Theorem~\ref{thm:weight_optimization_general}, local scaling from Lemma~\ref{lem:local_weights}, and unitary edge weights with fixed over-relaxation parameter $\alpha$. The ADMM parameters for the scaled algorithms are computed from Theorem~\ref{thm:optimal_f_alpha}. The step-sizes for the unscaled algorithms are empirically chosen.}
\label{fig:dqp}
\end{figure}

As a first example, we consider a distributed quadratic programming problem with $3$ agents, a decision vector $x\in\mathbf{R}^{4}$, and an objective function on the form $f(x) = \sum_{i\in\mathcal{V}} 1/2 x^\top Q_i x + q_i^\top x $ with
\small
\begin{equation*}
\begin{aligned}
Q_1 &= \begin{bmatrix}
    0.4236  & -0.0235 &  -0.0411 &   0.0023\\
   -0.0235  &  0.0113 &   0.0023 &  -0.0001\\
   -0.0411  &  0.0023 &   0.4713 &  -0.0262\\
    0.0023  & -0.0001 &  -0.0262 &   0.0115\\
\end{bmatrix}\\
Q_2 &= \begin{bmatrix}
    0.8417  & -0.1325  & -0.0827  &  0.0132\\
   -0.1325  &  0.0311  &  0.0132  & -0.0021\\
   -0.0827  &  0.0132  &  0.9376  & -0.1477\\
    0.0132  & -0.0021  & -0.1477  &  0.0335
\end{bmatrix}\\
Q_3 &=
\begin{bmatrix}
    0.0122  &  0.0308 &  -0.0002  & -0.0031\\
    0.0308  &  0.4343 &  -0.0031  & -0.0422\\
   -0.0002  & -0.0031 &   0.0125  &  0.0344\\
   -0.0031  & -0.0422 &   0.0344  &  0.4833
\end{bmatrix}\\
q_1&=q_2=0,\, q_3^\top = \begin{bmatrix}-0.1258 & 0.0087 & 0.0092 & -0.1398\end{bmatrix}.
\end{aligned}
\end{equation*}
\normalsize
The communication graph is a line graph where node $2$ is connected to nodes $1$ and $3$. The distributed optimization problem is formulated using edge variables and solved by executing the resulting ADMM iterations. The convergence behavior of the iterates for different choices of scalings and algorithm parameters are presented in Figure~\ref{fig:dqp}.

The optimal constraint scaling matrix and ADMM parameters are computed using Algorithm~\ref{alg:optimal_scaling_general}, resulting in $\rho = 1$ and $\alpha = 1.33$.
In the ``local'' algorithm, nodes determined constraint scalings in a distributed manner in accordance to Lemma~\ref{lem:local_weights}, while the optimal parameters computed using Theorem~\ref{thm:optimal_f_alpha} are $\rho = 1.44$ and $\alpha = 1.55$.
The remaining iterations correspond to ADMM algorithms with unitary edge weights, fixed relaxation parameter $\alpha$, and manually optimized step-size $\rho$. The parameter $\alpha$ is fixed at $1.0$, $1.5$, and $1.8$, while the corresponding $\rho$ is chosen as the empirical best.

Figure~\ref{fig:dqp} shows that the manually tuned ADMM algorithm exhibits worse performance than the optimally and locally scaled algorithms. Here, the best parameters for the scaled versions are computed systematically using the results derived earlier, while the best parameters for the unscaled algorithms are computed through exhaustive search.

\subsection{Distributed consensus}
In this section we apply our methodology to derive optimally scaled ADMM iterations for a particular problem instance usually referred to as \emph{average consensus}. The problem amounts to devising a distributed algorithm that ensures that all agents $i\in {\mathcal V}$ in a network reach agreement on the network-wide average of scalars $q_i$ held by the individual agents. This problem can be formulated as a particular case of~\eqref{eq:DQP_consensus} where $x\in \mathbf{R}$ and $f(x)=\sum_{i\in\mathcal{V}} 1/2(x-q_i)^2 = \sum_{i\in\mathcal{V}} 1/2x^2 -q_i x + 1/2q_i^2$. We consider edge-variable and node-variable formulations and compare the performance of the corresponding ADMM iterates with the relevant state-of-the-art algorithms.  As performance indicator, we use the convergence factors computed as the second largest eigenvalue of the linear fixed point iterations associated with each method. We generated communication graphs from the Erd\H{o}s-R\'enyi and the Random Geometric Graph (RGG) families (see,~\eg, \cite{Pen:2003}). Having generated $\vert \mathcal{V}\vert$ number of nodes, in Erd\H{o}s-R\'enyi graphs we connected each pair of nodes with probability $p=(1+\epsilon){\log(\vert\mathcal{V}\vert)}/{\vert\mathcal{V}\vert}$ where $\epsilon \in (0,1)$. In RGG,  $\vert\mathcal{V}\vert$ nodes were randomly deployed in the unit square and an edge was introduced between each pair of nodes whose inter-distance is at most $2\log(\vert\mathcal{V}\vert)/\vert\mathcal{V}\vert$; this guarantees that the graph is connected with high probability~\cite{kumar:00}.

\begin{figure}[p!tb]
  \centering
  \subfigure[RGG edge variable]
  {\includegraphics[width=0.4\hsize]{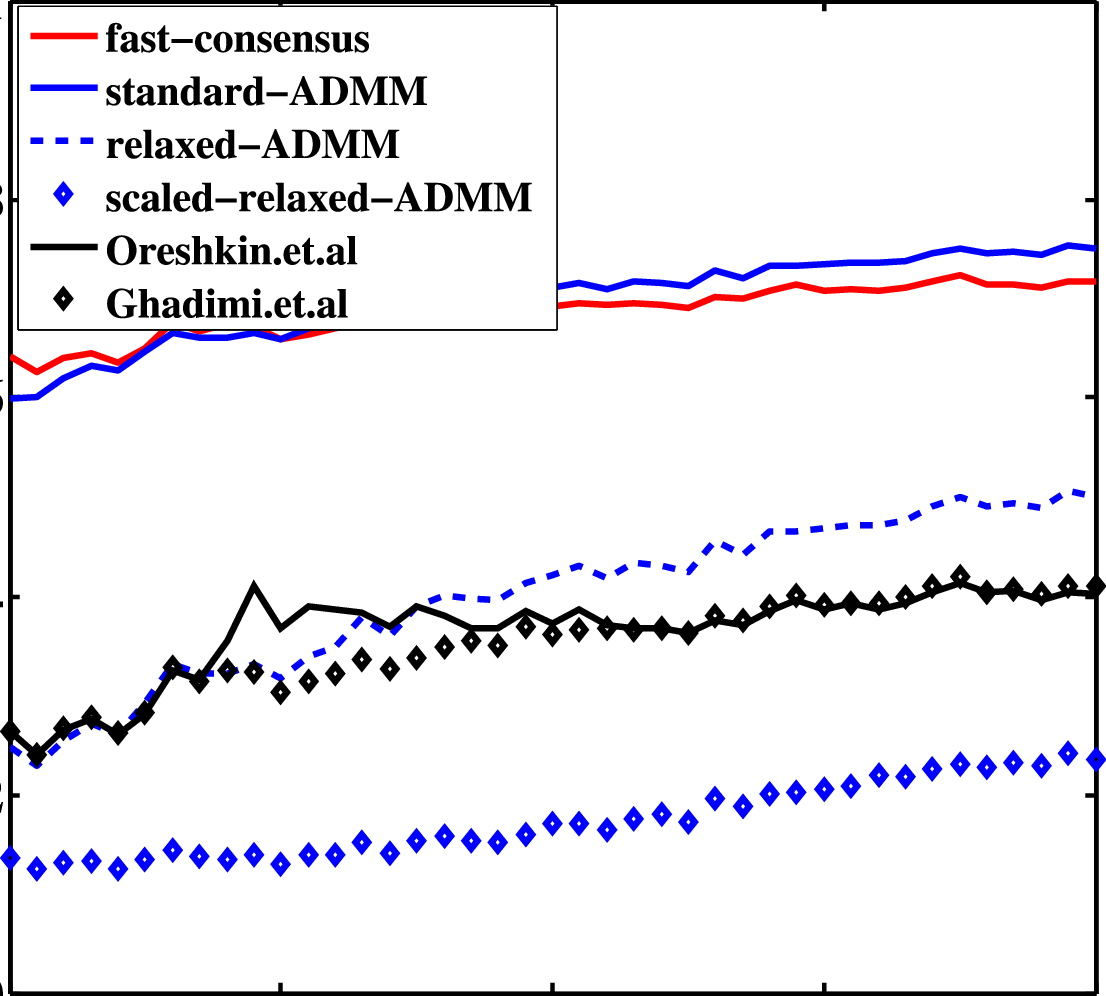}
  \label{fig:cons_0_1}}
  \subfigure[RGG node variable]
  {\includegraphics[width=0.4\hsize]{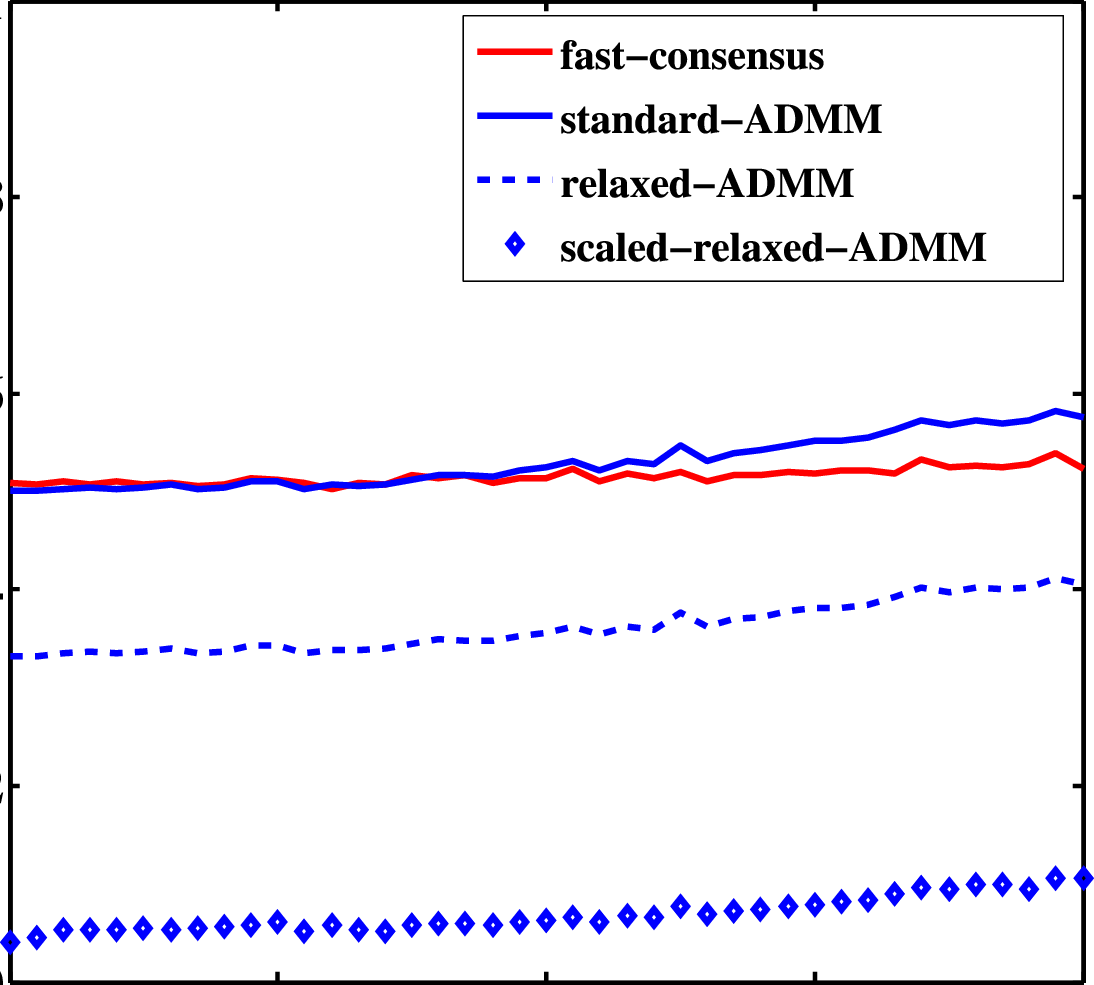}
  \label{fig:cons_0_2}}
  \subfigure[$\epsilon=0.2$ edge variable]
  {\includegraphics[width=0.4\hsize]{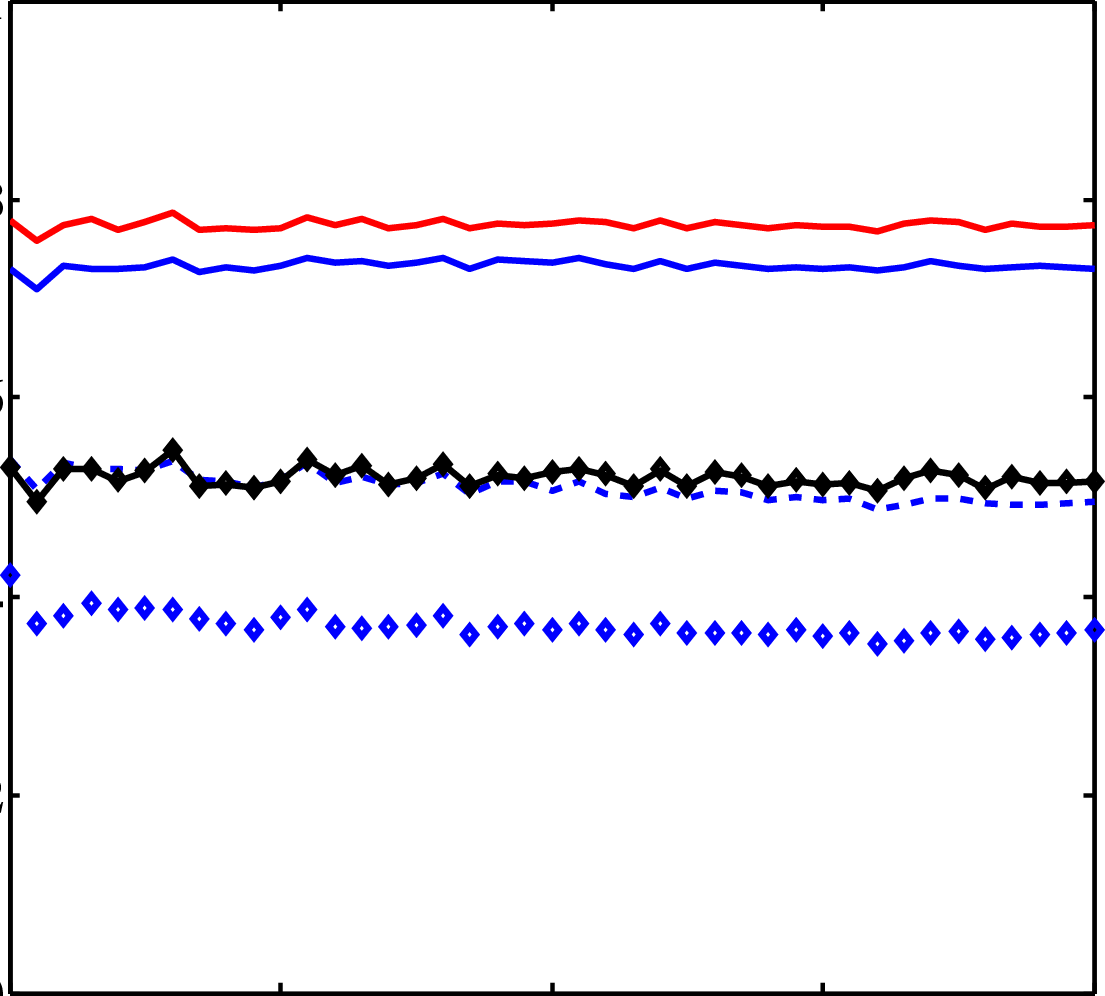}
  \label{fig:cons_1_1}}
  \subfigure[$\epsilon=0.2$ node variable]
  {\includegraphics[width=0.4\hsize]{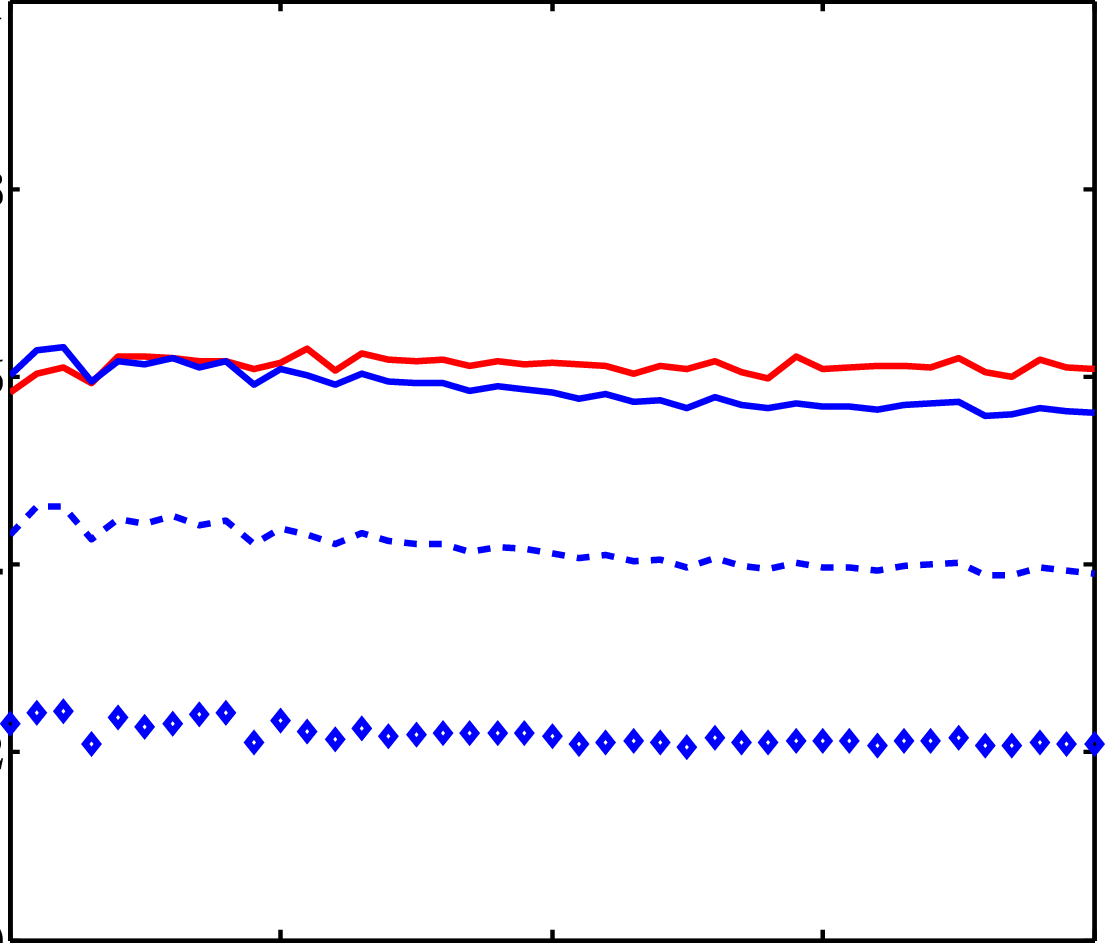}
  \label{fig:cons_1_2}}
  \subfigure[$\epsilon=0.8$ edge variable]
  {\includegraphics[width=0.4\hsize]{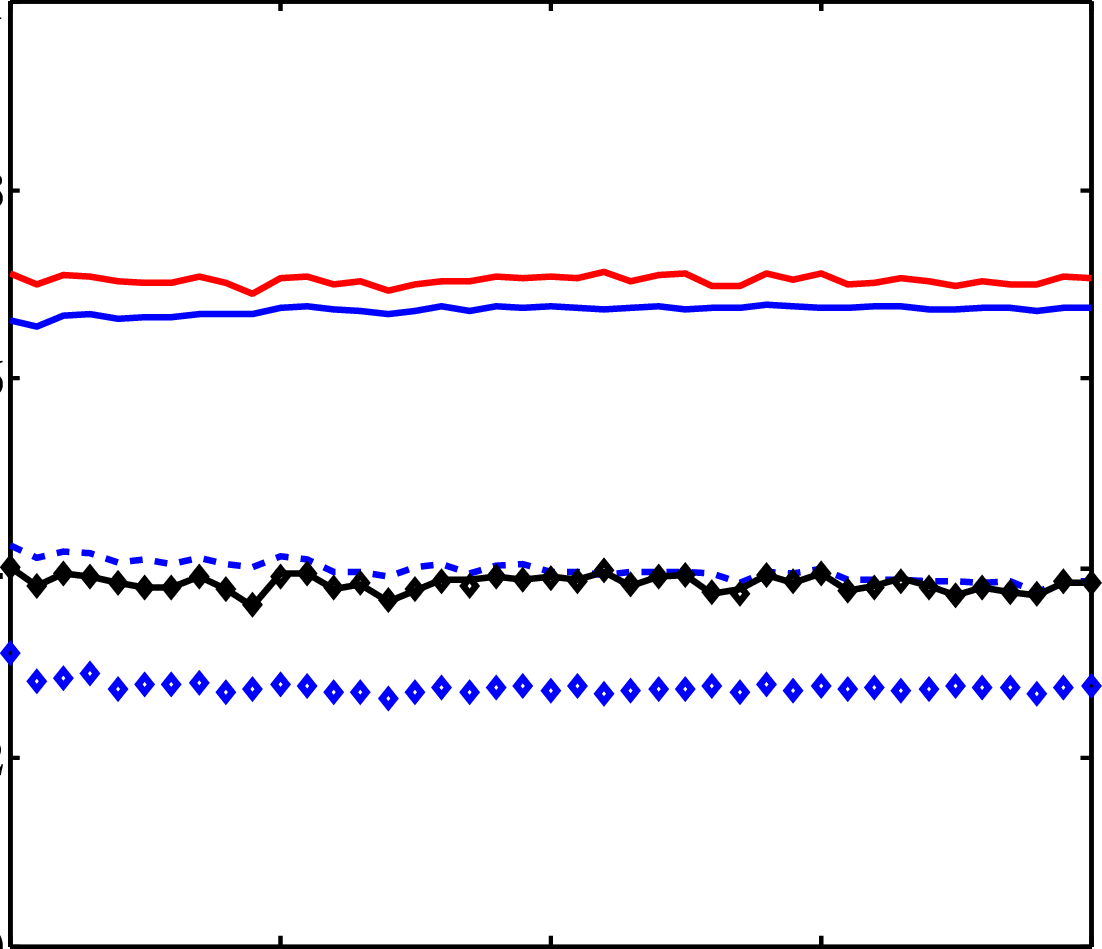}
  \label{fig:cons_2_1}}
  \subfigure[$\epsilon=0.8$ node variable]
  {\includegraphics[width=0.4\hsize]{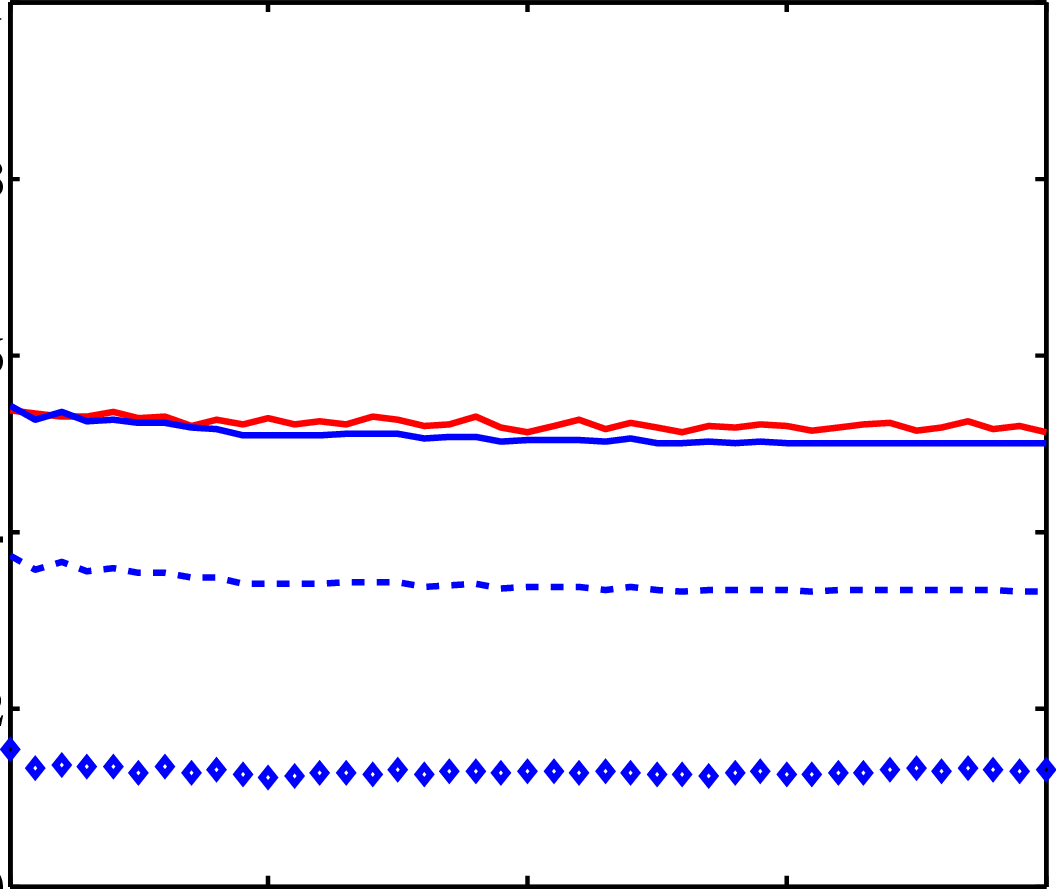}
  \label{fig:cons_2_2}}

\caption{Performance comparison of the proposed optimal scaling for the ADMM algorithm with state-of-the-art algorithms fast-consensus~\cite{Italian}, Oreshkin et.al~\cite{rabbat010} and Ghadimi et.al~\cite{ghadimi:13} The network of size $n=[10,50]$ is randomly generated by RGG (random geometric graphs) and Erd\H{o}s-R\'enyi graphs with low and high densities $\epsilon = \{0.2,0.8\}$.}
\label{fig:cons_1}
\end{figure}

Figure~\ref{fig:cons_1} presents Monte Carlo simulations of the convergence factors versus the number of nodes $\vert\mathcal{V}\vert\in [10, 50]$. Each data point is the average convergence factor in $60$ instances of randomly generated graphs with the same number of nodes.  In our simulations, we consider both edge-variable and node-variable formulations. For both formulations, we consider three versions of the ADMM algorithm with our parameter settings: the standard one (with step-size given in Corollary~\ref{cor:optimal_parameters_alpha_1}), an over-relaxed version with parameters in Theorem~\ref{thm:optimal_f_alpha}, and the scaled-relaxed-ADMM that uses weight optimization in addition to the optimal parameters in Theorem~\ref{thm:optimal_f_alpha}.

In the edge-variable scenario, we compare the ADMM iterates to three other algorithms: fast-consensus~\cite{Italian} from the ADMM literature and two state-of-the-art algorithms from the literature on the accelerated consensus: Oreshkin et al.~\cite{rabbat010} and Ghadimi et al.\cite{ghadimi:13}. In these algorithms, a two-tap memory mechanism is implemented so that the values of two last iterates are taken into account when computing the next. All the competitors employ the best weight scheme known for the respective method. For Ghadimi et al., the optimal weight is given in~\cite{ghadimi:13} while fast-consensus and Oreshkin et al. use the optimal weights in~\cite{XiaoBoyd04}. The scaled-relaxed-ADMM method employs the weight heuristic presented in  Lemma~\ref{lem:weight_optimization_edge}. Figures~\ref{fig:cons_0_1}, \ref{fig:cons_1_1} and~\ref{fig:cons_2_1} show a significant improvement of our design rules compared to the alternatives for both RGG and Erd\H{o}s-R\'enyi graphs in sparse ($\epsilon = 0.2$) and dense ($\epsilon = 0.8$) topologies. We observe that in all cases, the convergence factor decreases with increasing network size on RGG, while it stays almost constant on Erd\H{o}s-R\'enyi graphs.

For the node-variable formulation, we compare the three variants of our ADMM algorithm to the fast-consensus~\cite{Italian} algorithm. The reason why we exclude two other methods from the comparison is that they do not (yet) exist for the node-variable formulation. By comparing their explicit $x$-updates in~\eqref{eq:ADMM_edge_i} and~\eqref{eq:ADMM_node_i}, it is apparent that while each iterate of the consensus algorithms based on edge-variable formulation requires a single message exchange within the neighborhood of each node, the node-variable based algorithms require at least twice the number of message exchanges per round.
 In the scaled-relaxed-ADMM method, we first minimize the second largest generalized eigenvalue of $(2AD^{-1}A-D,D)$ using the quasi-convex program~\eqref{eqn:weight_optimization_general}. Let $A^\star$ and $D^\star$ be the adjacency and the degree matrices associated with the optimal solution of \eqref{eqn:weight_optimization_general}.  After extensive simulations it is observed that formulating the fixed point equation~\eqref{eq:ADMM_eq_matrix_x} as
\begin{equation}\label{eq:Consensus_node_linearMatrix_mod}
\begin{aligned}
 M_{11} &= \alpha\rho(Q+\rho D^\star)^{-1}(A^\star {D^\star}^{-1}A^\star)+ I,
 \iftoggle{draft}{%
\;   M_{12} =   \alpha\rho(Q+\rho D^\star)^{-1}, \; M_{21} = - \dfrac{\alpha}{2}(A^\star {D^\star}^{-1}A^\star+D^\star),
}{%
\\  M_{12} &=   \alpha\rho(Q+\rho D^\star)^{-1}, \quad M_{21} = - \dfrac{\alpha}{2}(A^\star {D^\star}^{-1}A^\star+D^\star),
}
\end{aligned}
\end{equation}
instead of using \eqref{eq:Consensus_node_linearMatrix}, often significantly improves the convergence factor of the ADMM algorithm for the node-variable formulation.  Note that this reformulation leads to $\{\lambda_i\}$ (in Theorem~\ref{thm:optimal_f_alpha}) being the generalized eigenvalues of $(AD^{-1}A,\;D)$. These eigenvalues have several nice properties,
\eg, they are positive and satisfy {Case I}, for which we presented the optimal ADMM parameters $(\alpha^\star, \rho^\star)$ in Theorem~\ref{thm:optimal_f_alpha}. The algorithm formulated by \eqref{eq:Consensus_node_linearMatrix_mod} corresponds to running the ADMM algorithm over a network with possible self loops  with the adjacency matrix $\tilde{A}=\tilde{A}^\top$ such that
\begin{align*}
A^\star {D^\star}^{-1  } A^\star = \dfrac{1}{2} (\tilde{A} {D^\star}^{-1} \tilde{A}+D^\star),\;
 \mbox{diag}(\tilde{A}\textbf{1}_n) = D^\star.
 \end{align*}
Figures~\ref{fig:cons_0_2}, \ref{fig:cons_1_2} and \ref{fig:cons_2_2} illustrate the performance benefits of employing optimal parameter settings developed in this paper compared to the alternative fast-consensus for different random topologies.

\section{Conclusions and future work}\label{sec:conclusion}
We considered the problem of optimal parameter selection and scaling of the ADMM method for distributed quadratic programming. Distributed unconstrained quadratic problems were cast as equality-constrained quadratic problems, to which the scaled ADMM method is applied. For this class of problems, the network-constrained scaling corresponds to the usual step-size constant, the relaxation parameter, and the edge weights of the communication graph. For connected communication graph, analytical expressions for the optimal step-size, relaxation parameter, and the resulting convergence factor were derived in terms of the spectral properties of the graph. Supposing the optimal step-size and relaxation parameter are chosen, the convergence factor is further minimized by optimally choosing the edge weights. Our results were illustrated in numerical examples and significant performance improvements over state-of-the-art techniques were demonstrated. As a future work, we plan to extend the results to a broader class of distributed quadratic problems.
\section*{Acknowledgment}
\addcontentsline{toc}{section}{Acknowledgment}
The authors would like to thank Michael Rabbat and Themistoklis Charalambous for valuable discussions and helpful comments to this paper.
\bibliography{admmbib}

\begin{thebibliography}{10}
\providecommand{\url}[1]{#1}
\csname url@samestyle\endcsname
\providecommand{\newblock}{\relax}
\providecommand{\bibinfo}[2]{#2}
\providecommand{\BIBentrySTDinterwordspacing}{\spaceskip=0pt\relax}
\providecommand{\BIBentryALTinterwordstretchfactor}{4}
\providecommand{\BIBentryALTinterwordspacing}{\spaceskip=\fontdimen2\font plus
\BIBentryALTinterwordstretchfactor\fontdimen3\font minus
  \fontdimen4\font\relax}
\providecommand{\BIBforeignlanguage}[2]{{%
\expandafter\ifx\csname l@#1\endcsname\relax
\typeout{** WARNING: IEEEtran.bst: No hyphenation pattern has been}%
\typeout{** loaded for the language `#1'. Using the pattern for}%
\typeout{** the default language instead.}%
\else
\language=\csname l@#1\endcsname
\fi
#2}}
\providecommand{\BIBdecl}{\relax}
\BIBdecl

\bibitem{nedic10}
A.~Nedic, A.~Ozdaglar, and P.~Parrilo, ``Constrained consensus and optimization
  in multi-agent networks,'' \emph{Automatic Control, IEEE Transactions on},
  vol.~55, no.~4, pp. 922--938, Apr. 2010.

\bibitem{Italian}
T.~Erseghe, D.~Zennaro, E.~Dall'Anese, and L.~Vangelista, ``Fast consensus by
  the alternating direction multipliers method,'' \emph{Signal Processing, IEEE
  Transactions on}, vol.~59, pp. 5523--5537, 2011.

\bibitem{accelerated_mpc13}
P.~Giselsson, M.~D. Doan, T.~Keviczky, B.~D. Schutter, and A.~Rantzer,
  ``Accelerated gradient methods and dual decomposition in distributed model
  predictive control,'' \emph{Automatica}, vol.~49, no.~3, pp. 829--833, 2013.

\bibitem{farhad12}
F.~Farokhi, I.~Shames, and K.~H. Johansson, ``Distributed {MPC} via dual
  decomposition and alternative direction method of multipliers,'' in
  \emph{Distributed Model Predictive Control Made Easy}, ser. Intelligent
  Systems, Control and Automation: Science and Engineering, J.~M. Maestre and
  R.~R. Negenborn, Eds.\hskip 1em plus 0.5em minus 0.4em\relax Springer, 2013,
  vol.~69.

\bibitem{Falcao1995}
D.~Falcao, F.~Wu, and L.~Murphy, ``Parallel and distributed state estimation,''
  \emph{Power Systems, IEEE Transactions on}, vol.~10, no.~2, pp. 724--730, May
  1995.

\bibitem{Boyd11}
S.~Boyd, N.~Parikh, E.~Chu, B.~Peleato, and J.~Eckstein, ``Distributed
  optimization and statistical learning via the alternating direction method of
  multipliers,'' \emph{Foundations and Trends in Machine Learning}, vol. 3
  Issue: 1, pp. 1--122, 2011.

\bibitem{tyler12}
C.~Conte, T.~Summers, M.~Zeilinger, M.~Morari, and C.~Jones, ``Computational
  aspects of distributed optimization in model predictive control,'' in
  \emph{Decision and Control (CDC), 2012 IEEE 51st Annual Conference on}, 2012.

\bibitem{marriette12}
M.~Annergren, A.~Hansson, and B.~Wahlberg, ``An {ADMM} algorithm for solving $
  \ell_1 $ regularized {MPC},'' in \emph{Decision and Control (CDC), 2012 IEEE
  51st Annual Conference on}, 2012.

\bibitem{Mota12}
J.~Mota, J.~Xavier, P.~Aguiar, and M.~Puschel, ``Distributed admm for model
  predictive control and congestion control,'' in \emph{Decision and Control
  (CDC), 2012 IEEE 51st Annual Conference on}, 2012.

\bibitem{luo12}
Z.~Luo, ``On the linear convergence of the alternating direction method of
  multipliers,'' \emph{ArXiv e-prints}, 2012.

\bibitem{boley12}
D.~Boley, ``Local linear convergence of the alternating direction method of
  multipliers on quadratic or linear programs,'' \emph{SIAM Journal on
  Optimization}, vol.~23, pp. 2183--2207, 2013.

\bibitem{deng12}
W.~Deng and W.~Yin, ``On the global and linear convergence of the generalized
  alternating direction method of multipliers,'' Rice University CAAM Technical
  Report ,TR12-14, 2012., Tech. Rep., 2012.

\bibitem{GTS:13}
E.~Ghadimi, A.~Teixeira, I.~Shames, and M.~Johansson, ``Optimal parameter
  selection for the alternating direction method of multipliers ({ADMM}):
  quadratic problems,'' \emph{IEEE Transactions on Automatic Control}, 2014, to
  appear.

\bibitem{Gomez2011_DSE}
A.~G\'{o}mez-Exp\'{o}sito, A.~de~la Villa~Ja\'{e}n, C.~G\'{o}mez-Quiles,
  P.~Rousseaux, and T.~V. Cutsem, ``A taxonomy of multi-area state estimation
  methods,'' \emph{Electric Power Systems Research}, vol.~81, no.~4, pp.
  1060--1069, 2011.

\bibitem{teixeira_admm_2013}
A.~Teixeira, E.~Ghadimi, I.~Shames, H.~Sandberg, and M.~Johansson, ``Optimal
  scaling of the admm algorithm for distributed quadratic programming,'' in
  \emph{Proceedings of the {IEEE} 52nd Conference on Decision and Control},
  Dec. 2013, pp. 6868--6873.

\bibitem{ghadimi2014admm}
E.~Ghadimi, A.~Teixeira, M.~Rabbat, and M.~Johansson, ``The {ADMM} algorithm
  for distributed averaging: Convergence rates and optimal parameter
  selection,'' in \emph{Proceedings of the 48th Asilomar Conference on Signals,
  Systems and Computers}, 2014, to appear.

\bibitem{NEO:09}
A.~Nedic and A.~Ozdaglar, ``Distributed subgradient methods for multi-agent
  optimization,'' \emph{Automatic Control, IEEE Transactions on}, vol.~54,
  no.~1, pp. 48--61, Jan 2009.

\bibitem{Jury_stability_test}
E.~Jury, \emph{Theory and Application of the z-Transform Method}.\hskip 1em
  plus 0.5em minus 0.4em\relax Huntington, New York: Krieger Publishing
  Company, 1974.

\bibitem{chung1997spectral}
F.~R. Chung, \emph{Spectral graph theory}.\hskip 1em plus 0.5em minus
  0.4em\relax American Mathematical Soc., 1997, vol.~92.

\bibitem{butler2008eigenvalues}
S.~K. Butler, \emph{Eigenvalues and structures of graphs}.\hskip 1em plus 0.5em
  minus 0.4em\relax University of California, San Diego, ProQuest, UMI
  Dissertations Publishing, 2008.

\bibitem{Pen:2003}
M.~Penros, \emph{Random Geometric Graphs}.\hskip 1em plus 0.5em minus
  0.4em\relax Oxford Studies in Probability, 2003.

\bibitem{kumar:00}
P.~Gupta and P.~Kumar, ``The capacity of wireless networks,'' \emph{Information
  Theory, IEEE Transactions on}, vol.~46, no.~2, pp. 388--404, Mar 2000.

\bibitem{rabbat010}
B.~Oreshkin, M.~Coates, and M.~Rabbat, ``Optimization and analysis of
  distributed averaging with short node memory,'' \emph{Signal Processing, IEEE
  Transactions on}, vol. 58 Issue: 5, pp. 2850--2865, 2010.

\bibitem{ghadimi:13}
E.~Ghadimi, I.~Shames, and M.~Johansson, ``Multi-step gradient methods for
  networked optimization,'' \emph{Signal Processing, IEEE Transactions on},
  vol.~61, no.~21, pp. 5417--5429, Nov 2013.

\bibitem{XiaoBoyd04}
L.~Xiao and S.~Boyd, ``Fast linear iterations for distributed averaging,''
  \emph{Systems and Control Letters}, vol. 53 Issue: 1, pp. 65--78, 2004.

\end{thebibliography}

\appendix

\subsection{Proof of Lemma~\ref{lem:QP_change_of_variable}}\label{app:lem:QP_change_of_variable}
Rewrite~\eqref{eq:QP_problem} in terms of the variables $\tilde{x}=x-\hat{x}$ and $\tilde{z}=z-\hat{z}$:
  \begin{align*}
  \begin{array}[c]{ll}
	\underset{\tilde{x},\tilde{z}}{\mbox{minimize}} &\dfrac{1}{2}(\tilde{x}+\hat{x})^\top Q(\tilde{x}+\hat{x}) + q^{\top}(\tilde{x}+\hat{x})+c^\top (\tilde{z}+\hat{z}) \\
	\mbox{subject to} &
    RE(\tilde{x}+\hat{x}) + RF(\tilde{z}+\hat{z}) =R h.
	\end{array}
    \end{align*}
  Collecting the terms in the objective and noting that the feasible solution $(\hat{x},\hat{z})$ satisfies $RE\hat{x}+RF\hat{z}=Rh$, one can rewrite this problem as
  \begin{equation}\nonumber
    \begin{aligned}
     \begin{array}[c]{ll}
	\underset{\tilde{x},\tilde{z}}{\mbox{minimize}} &\dfrac{1}{2}\tilde{x}^\top Q\tilde{x} + (Q\hat{x}+q)^\top\tilde{x}+c^\top \tilde{z}+ d \\
	\mbox{subject to} &
    RE \tilde{x} + RF \tilde{z} = 0,
	\end{array}
    \end{aligned}
  \end{equation}
  where $d$ collects the constant terms. Since $d$ does not affect the minimizer, it can be removed and the problem is equivalent to~\eqref{eq:QP_problem_suitable} when $p= Q\hat{x}+q$.
\subsection{Proof of Theorem~\ref{thm:convergence_ADMM}}\label{app:thm:convergence_ADMM}
Consider the quadratic programming problem~\eqref{eq:QP_problem} with $\bar{E} = RE$, $\bar{F}=RF$ and $h=0$. Defining the feasibility subspace  as ${\mathcal{X} \triangleq \left\{x\in\mathbf{R}^n,\, z \in\mathbf{R}^m\vert\; \bar{E}x+ \bar{F}z = 0  \right\}}$, the dimension of $\mathcal{X}$ is given by $\mbox{dim}\left(\mathcal{X}\right) = \mbox{dim}\left(\Null{[\bar{E} \, \bar{F}]}\right)$. Observe that we have $\mbox{dim}\left(\Range{\bar{E}}\right) = n$ and $\mbox{dim}\left(\Range{\bar{F}}\right) = m$, since $\bar{E}\in\mathbf{R}^{r\times n}$ and $\bar{F}\in\mathbf{R}^{r\times m}$ have full column rank. Using the equalities $\mbox{dim}\left(\mathcal{R}([\bar{E} \, \bar{F}])\right) = \mbox{dim}\left(\Range{\bar{E}}\right) + \mbox{dim}\left(\Range{\bar{F}}\right) - \mbox{dim}\left(\Range{\bar{E}}\cap\Range{\bar{F}} \right)$ and $\mbox{dim}\left(\Null{[\bar{E} \, \bar{F}]}\right) + \mbox{dim}\left(\Range{[\bar{E} \, \bar{F}]}\right) = n+m$, we conclude that
$\mbox{dim}\left(\mathcal{X}\right) = \mbox{dim}\left(\Range{\bar{E}}\cap\Range{\bar{F}} \right) = s$.

Provided that~\eqref{eq:QP_problem} is feasible and under the assumption that there exists (non-trivial) non-zero tuple $(x,z)\in \mathcal{X}$,  we have $s\geq 1$. A necessary condition for the ADMM iterations to converge to a fixed-point $(x^\star,z^\star)$ is that $M$ (in the fixed-point iterates $\sigma^{k+1}=M\sigma^k$) has $\phi_{2n}=1$.
Moreover, when $\phi_{2n}=1$ the ADMM iterations converge to the $1$-eigenspace of $M$ defined as $\mbox{span}(V)$ with $MV = V$, where the dimension of $\mbox{span}(V)$ corresponds to the multiplicity of the $1$-eigenvalue.

Given a feasibility subspace $\mathcal{X}$ different problem parameters $Q, q$, and $c$ lead to different optimal solution points in $\mathcal{X}$. Therefore, for the fixed-point $(x^\star, z^\star)$ to be the optimal, the span of the fixed-points of $M$ must contain the whole feasibility subspace $\mathcal{X}$. That is,
the $1$-eigenvalue must have multiplicity $\mbox{dim}(\mbox{span}(V)) = \mbox{dim}\left(\mathcal{X}\right) = s$, \ie, $1=\phi_{2n} = \dots = \phi_{2n-s+1} > |\phi_{2n-s}|$.

Next we show that fixed-points of the ADMM iterations satisfy the optimality conditions of~\eqref{eq:QP_problem} in terms of the augmented Lagrangian. The fixed-point of the ADMM iterations~\eqref{eqn:admm_scaled_iterations} satisfy the system of equations
\begin{equation}\label{eqn:fixed_point_optimality_cond}
\begin{aligned}
\begin{bmatrix}
Q+\rho \bar{E}^\top \bar{E} & \rho\bar{E}^\top \bar{F} & \rho\bar{E}^\top \\
\alpha \bar{F}^\top \bar{E} & \alpha \bar{F}^\top\bar{F} & \bar{F}^\top\\
\bar{E} & \bar{F} & 0
\end{bmatrix}
\begin{bmatrix}
x^{\star}\\
z^{\star}\\
u^{\star}
\end{bmatrix}
=
\begin{bmatrix}
-q\\
- {c/\rho}\\
0
\end{bmatrix}.
\end{aligned}
\end{equation}

From Karush-Kuhn-Tucker optimality conditions of the augmented Lagrangian $L_\rho(x,z,u)=1/2 x^\top Q x + q^\top x + c^\top z+ \rho/2\Vert \bar{E}x+\bar{F}z\Vert^2 + \rho u^\top (\bar{E}x+\bar{F}z)$ it yields
\begin{equation*}
\begin{aligned}
\begin{bmatrix}
Q+\rho \bar{E}^\top \bar{E} & \rho\bar{E}^\top \bar{F} & \rho\bar{E}^\top \\
\rho \bar{F}^\top \bar{E} & \rho \bar{F}^\top\bar{F} & \rho\bar{F}^\top\\
\bar{E} & \bar{F} & 0
\end{bmatrix}
\begin{bmatrix}
x^{\star}\\
z^{\star}\\
u^{\star}
\end{bmatrix}
=
\begin{bmatrix}
-q\\
-c\\
0
\end{bmatrix},
\end{aligned}
\end{equation*}
which is equivalent to~\eqref{eqn:fixed_point_optimality_cond} by noting that $\bar{F}^\top \bar{E} x^\star+ \bar{F}^\top \bar{F}z^\star = \bar{F}^\top (\bar{E}x^\star+\bar{F}z^\star)=0$.

\subsection{Proof of Theorem~\ref{thm:M_eigenvalues}}\label{app:thm:M_eigenvalues}
To satisfy the eigenvalue equation $M [v_i^\top \; w_i^\top]^\top = \phi_i [v_i^\top \; w_i^\top]^\top$, $v_i$ and $w_i$ should satisfy
\begin{equation}\nonumber
  \begin{aligned}
    & \left( M_{11} +  \dfrac{1}{\phi_i-1+\alpha} M_{12}M_{21} -\phi_i I  \right)v_i=0,\\
    & w_i = \dfrac{1}{(\phi_i-1+\alpha)} M_{21} v_i.
  \end{aligned}
\end{equation}
When $\bar{E}^\top \bar{E} = \kappa Q$, we have
\begin{equation}\nonumber
\begin{aligned}
  &\left( M_{11} +  \dfrac{1}{\phi_i-1+\alpha} M_{12}M_{21} -\phi_i I  \right)v_i \\
  &=\alpha\beta(\bar{E}^\top \bar{E})^{-1}\bar{E}^\top\left(2\Pi_{\mathcal{R}(\bar{F})} -  I\right)\bar{E} v_i+v_i \\
   &  -  \dfrac{\alpha^2 \beta}{\phi_i-1+\alpha}(\bar{E}^\top \bar{E})^{-1} \bar{E}^\top \Pi_{\mathcal{R}(\bar{F})}  \bar{E} v_i - \phi_i v_i\\
   & = (\alpha \beta \lambda_i+1) v_i - \dfrac{\alpha^2\beta}{2}\dfrac{\lambda_i+1}{\phi_i-1+\alpha}v_i-\phi_i v_i \\
   &= \left(\alpha \beta \lambda_i+1-\phi_i - \dfrac{\alpha^2\beta(\lambda_i+1)}{2(\phi_i-1+\alpha)}\right)v_i,
\end{aligned}
\end{equation}
where the last steps follow from the generalized eigenvalue assumption. Thus, the eigenvalues of $M$ are given as the solution of
\begin{equation}\nonumber
  \phi_i^2 + (\alpha-\alpha\beta \lambda_i-2)\phi_i + \alpha \beta \lambda_i(1-\dfrac{\alpha}{2}) + \dfrac{1}{2}\alpha^2 \beta+1-\alpha=0.
\end{equation}
\subsection{Proof of Lemma~\ref{lem:bar_lambda}}\label{app:lem:bar_lambda}
Recall that a complex number $\lambda_i$ is a generalized eigenvalue of $\left(\bar{E}^\top (2\Pi_{\mathcal{R}(\bar{F})} - I) \bar{E},\;\bar{E}^\top \bar{E} \right)$ if
there exists a non-zero vector $\nu_i\in\mathbf{C}^n$ such that $ \left(\bar{E}^\top (2\Pi_{\mathcal{R}(\bar{F})} - I) \bar{E} - \lambda_i \bar{E}^\top \bar{E} \right) \nu_i = 0$. Since $\bar{E}$ has full column rank, $\bar{E}^\top \bar{E}$ is invertible and we observe that $\lambda_i$ is an eigenvalue of the symmetric matrix $  (\bar{E}^\top \bar{E})^{-1/2}\bar{E}^\top (2\Pi_{\mathcal{R}(\bar{F})} - I) \bar{E} (\bar{E}^\top \bar{E})^{-1/2}$. Since the latter is a real symmetric matrix, we conclude that the generalized eigenvalues and eigenvectors are real.

For the second part of the proof, note that the following bounds hold for a generalized eigenvalue $\lambda_i$
\begin{equation*}
\min_{\nu\in\mathbf{R}^n} \dfrac{2\nu^\top \bar{E}^\top \Pi_{\mathcal{R}(\bar{F})} \bar{E} \nu}{\nu^\top \bar{E}^\top \bar{E}\nu}-1  \leq \lambda_i \leq \max_{\nu\in\mathbf{R}^n} \dfrac{2\nu^\top \bar{E}^\top \Pi_{\mathcal{R}(\bar{F})} \bar{E} \nu}{\nu^\top \bar{E}^\top \bar{E}\nu}-1.
\end{equation*}
Since the projection matrix $\PrjR{\bar{F}}$ only takes $0$ and $1$ eigenvalues we have $0\leq 2 \nu^\top \bar{E}^\top \PrjR{\bar{F}} \bar{E} \nu \leq 2\nu^\top \bar{E}^\top \bar{E} \nu$ which shows that $\lambda_i\in [-1,1]$.

\subsection{Proof of Lemma~\ref{lem:bar_lambda_ones}}\label{app:lem:bar_lambda_ones}
Let $V_\mathcal{X}\in\mathbf{R}^{(n+m)\times s}$ be a matrix whose columns are a basis for the feasibility subspace $\mathcal{X}$ and partition this matrix as $V_\mathcal{X} = [V_x^\top\, V_z^\top ]^\top$. We first show that the generalized eigenvectors associated with the unit generalized eigenvalues $\lambda_i=1$ are in $V_x$.

Given the partitioning of $V_\mathcal{X}$ we have that $\bar{E}V_x + \bar{F}V_z = 0$ and $\bar{E}\nu \in \Range{\bar{F}}$ for $\nu\in V_x$. Hence we have $\PrjR{\bar{F}} \bar{E} \nu = \bar{E} \nu$, yielding $(\nu^\top (\bar{E}^\top (2\PrjR{\bar{F}}-I)\bar{E})\nu)/(\nu^\top(\bar{E}^\top \bar{E})\nu)=1$. Moreover, as $1$ is the upper bound for ${\lambda_i}$ according to Lemma~\ref{lem:bar_lambda}, we conclude that $\lambda_n = 1$ is a generalized eigenvalue associated with the eigenvector $\nu$.
Next we derive the rank of $V_x$, which corresponds to the multiplicity of the unit generalized eigenvalue. Recall from the proof of Theorem~\ref{thm:convergence_ADMM} that the feasibility subspace $\mathcal{X}$ has $\mbox{dim}(\mathcal{X})=\mbox{dim}(\mathcal{R}(\bar{E}) \cap \mathcal{R}(\bar{F}))=s\geq 1$.
Given that $\bar{F}$ has full column rank, using the equation $\bar{E}V_x + \bar{F}V_z = 0$ we have that $V_z = -\bar{F}^\dagger \bar{E} V_x$. Hence, we conclude that $\mbox{rank}(V_\mathcal{X}) = \mbox{rank}(V_x) = s$ and that there exist $s$ generalized eigenvalues equal to $1$.

\subsection{Proof of Lemma~\ref{lem:admm_stability}}\label{app:lem:admm_stability}
Recall from Lemma~\ref{lem:bar_lambda_ones} that for a feasible problem of the form~\eqref{eq:QP_problem} we have $\lambda_i = 1$ for $i\geq n-s+1$. From~\eqref{eq:eigenvalues_M_phi} it follows that each $\lambda_i = 1$ results in two eigenvalues $\phi=1$  and $ \phi=1 - \alpha( 1-\beta ) $. Thus we conclude that $M$ has at least $s$ eigenvalues equal to $1$. Moreover, since $\beta\in(0,1)$ and $\alpha\in(0,2]$, we observe that $|1 - \alpha( 1-\beta )|<1$.
Next we consider $i < n-s+1$ and show that the resulting eigenvalues of $M$ are inside the unit circle for all $\beta\in(0,1)$ and $\alpha\in(0,2]$ using the necessary and sufficient conditions from Proposition~\ref{prop:Jury}.

The first condition of Proposition~\ref{prop:Jury} can be rewritten as $a_0+a_1+a_2=1/2\alpha^2\beta(1-\lambda_i) > 0$, which holds for $\lambda_i\in[-1, 1)$.
Having $\alpha > 0$ and $\lambda_i < 1$, the condition $a_2>a_0$ can be rewritten as
$\alpha < \left(2(1-\beta\lambda_i)\right) / \left(\beta(1-\lambda_i)\right)$. For $\beta<1$, that the right hand side term is greater than $2$, from which we conclude that the second condition is satisfied. It remains to show  $a_0 - a_1 + 1>0$. Replacing the terms on the left-hand-side, they form a convex quadratic polynomial on $\alpha$, \ie, $D(\alpha) = \dfrac{1}{2}\alpha^2 \beta(1-\lambda_i)  + 2\alpha ( \beta{\lambda_i}  -1) + 4$. The value of $\alpha$ minimizing $D(\alpha)$ is $\alpha = \left(2(1-\beta\lambda_i)\right) / \left(\beta(1-\lambda_i)\right)$, which was shown to be greater than $2$ when addressing the second condition. Since $D(2) = 2\beta(1+\lambda) >0$, we conclude $D(\alpha) > 0$ for all $\alpha \leq 2$ and the third condition holds.

\subsection{Proof of Theorem~\ref{thm:reduced_max}}\label{app:thm:reduced_max}
The magnitude of $\phi_{2n-s}$ can be characterized with Jury's stability test as follows. Consider the non-unit generalized eigenvalues $ \{\lambda_i\}_{i\leq n-s}$ and let $\phi_i = r\tilde{\phi_i}$ for $r\geq 0$. Substituting $\phi_i$ in the eigenvalue polynomials~\eqref{eq:eigenvalues_M_phi} yields  $r^2\tilde{\phi_i}^2 + ra_1(\lambda_i) \tilde{\phi_i} + a_0(\lambda_i)=0$.
Therefore, having the roots of these polynomials inside the unit circle
is equivalent to having $\vert \phi_i \vert < r$. From the stability of ADMM iterates (see Lemma~\ref{lem:admm_stability}) it follows that it is always possible to find $r<1$. Using the necessary and sufficient conditions from Proposition~\ref{prop:Jury}, $\vert\phi_{2n-s}\vert$ is obtained as
\begin{equation}\label{eq:phi_Jury}
%
\begin{aligned}
\underset{r\geq 0}{\mbox{minimize}} &\quad r &\\
\mbox{subject to}&\quad  a_0(\lambda_i) + r a_1(\lambda_i) + r^2 \geq 0&\\
& \quad  r^2 \geq a_0(\lambda_i)&\qquad \forall\, i\leq n-s\\
& \quad a_0(\lambda_i) - r a_1(\lambda_i) + r^2 \geq 0&\\
& \quad r \geq \vert1 - \alpha(1-\beta) \vert. &
\end{aligned}
\end{equation}

Next we remove redundant constraints from~\eqref{eq:phi_Jury}. Considering the first constraint, we aim at finding $\lambda\in\{\lambda_i\}_{i\leq n-s}$ such that $a_0(\lambda_i) + r a_1(\lambda_i) + r^2 \geq  a_0(\lambda) + r a_1(\lambda) + r^2$ for all $i\leq n-s$. Observing that the former inequality can be rewritten as $\alpha\beta(\lambda - \lambda_i)(1-\frac{\alpha}{2}-r)\leq 0$, we conclude that $\lambda = \lambda_{n-s}$ if $1-\frac{\alpha}{2} \leq r$ and $\lambda = \lambda_1$ otherwise. Hence the constraints $a_0(\lambda_i) + r a_1(\lambda_i) + r^2 \geq 0$ for $1<i<n-s$ are redundant.
As for the second condition, note that $a_0(\lambda_{n-s}) - a_0(\lambda_i) = \alpha\beta(\lambda_{n-s} - \lambda_i)(1-\frac{\alpha}{2}) \geq 0$ for all $i\leq n-s$, since $\alpha\in(0,\,2]$. Consequently, the constraints $r^2 \geq  a_0(\lambda_{i})$ for $i< n-s$ can be removed.
Regarding the third constraint, we aim at finding $\lambda\in\{\lambda_i\}_{i\leq n-s}$ such that $a_0(\lambda_i) - r a_1(\lambda_i) + r^2 \geq  a_0(\lambda) - r a_1(\lambda) + r^2$ for all $i\leq n-s$. Since the previous inequality can be rewritten as $\alpha\beta(\lambda - \lambda_i)(1-\frac{\alpha}{2}+r)\leq 0$, which holds for $\lambda = \lambda_{1}$, we conclude that the constraints for $1<i\leq n-s$ are redundant.
Removing the redundant constraints, the optimization problem~\eqref{eq:phi_Jury} can be rewritten as
\begin{equation}\label{eq:Jury:phi:minimization}
\begin{array}{ll}
\underset{r\geq 0,\{s_i\}}{\mbox{minimize}} &\quad r \\
\mbox{subject to}& \quad  a_0(\lambda_{n-s}) + r a_1(\lambda_{n-s}) + r^2 - s_1 = 0\\
&\quad  a_0(\lambda_{1}) + r a_1(\lambda_{1}) + r^2 - s_2 = 0\\
& \quad  r^2 - a_0(\lambda_{n-s}) -s_3 = 0  \\
& \quad a_0(\lambda_1) - r a_1(\lambda_1) + r^2 -s_4 =0\\
& \quad r - \vert1 - \alpha(1-\beta) \vert -s_5 =0 \\
& \quad s_i \geq 0 \quad \forall i\leq 5,
\end{array}
\end{equation}
where $\{s_i\}$ are slack variables. Subtracting the fourth equation from the second
we obtain the following equivalent problem
\begin{equation}\label{eq:phi_Jury_2}
\begin{array}{ll}
\underset{\{s_i\}}{\mbox{minimize}} &\quad \underset{i}{\max}\{r_i\} \\
\mbox{subject to}&\quad
s_i \geq 0, \quad \forall i\leq 5\\
& \quad r_i \geq 0, \quad \forall i\leq 7\\
& \quad a_0(\lambda_{n-s}) +s_3  \geq 0\\
& \quad \beta^2\lambda_{n-s}^2 - 2\beta + 1 + s_1 \geq 0\\
& \quad \beta^2\lambda_{1}^2 - 2\beta + 1 + s_4 \geq 0,
\end{array}
\end{equation}
where
\begin{equation*}
\begin{aligned}
r_1 &= 1 - \frac{\alpha}{2}  + \frac{\alpha}{2}\beta\lambda_{n-s} + \frac{\alpha}{2}\sqrt{\beta^2\lambda_{n-s}^2 - 2\beta + 1 + s_1  }\\
r_2 &= \dfrac{s_2-s_4}{2a_1(\lambda_1)}\\
 r_3 &= \sqrt{ a_0(\lambda_{n-s}) +s_3 }  \\
r_4 &= -1 + \frac{\alpha}{2}  - \frac{\alpha}{2}\beta\lambda_{1} + \frac{\alpha}{2}\sqrt{\beta^2\lambda_{1}^2 - 2\beta + 1 + s_4  }\\
r_5 &= \vert1 - \alpha(1-\beta) \vert +s_5 \\
 r_6 &= 1 - \frac{\alpha}{2}  + \frac{\alpha}{2}\beta\lambda_{n-s} - \frac{\alpha}{2}\sqrt{\beta^2\lambda_{n-s}^2 - 2\beta + 1 + s_1  }\\
r_7 &= -1 + \frac{\alpha}{2}  - \frac{\alpha}{2}\beta\lambda_{1} - \frac{\alpha}{2}\sqrt{\beta^2\lambda_{1}^2 - 2\beta + 1 + s_4  }.
\end{aligned}
\end{equation*}
In the above equation, $\{r_1,\, r_6 \}$, $r_2$, $r_3$, $\{r_4,\, r_7 \}$, and $r_5$ are solutions to the first, second, third, forth and fifth equality constraints in~\eqref{eq:Jury:phi:minimization}, respectively.
The last three inequalities impose that $r_1$, $r_3$, $r_4$, $r_6$, and $r_7$ are real values. Moreover, the last two constraints ensure that the inequalities $r_1\geq r_6$ and $r_4\geq r_7$ hold. Performing the minimization of each $r_i$ with respect to the corresponding slack variable $s_i$ we obtain
$\vert \phi_{2n-s} \vert = \max\{r^\star_1,\, r^\star_3,\, r^\star_4 ,\, r^\star_5\}$
where $r_i^\star$ are computed as in~\eqref{eq:phi_Jury_2} with
\begin{equation*}
\begin{aligned}
s_1^\star &= \max\{0,\; - (\beta^2\lambda_{n-s}^2 -2\beta +1  ) \},  \\
s_2^\star &=s_4^\star=  \max\{0,\; - (\beta^2\lambda_{1}^2 -2\beta +1  ) \}, \\
s_3^\star &=  \max\{0,\; - a_0(\lambda_{n-s})  \}, \quad
s_5^\star =  0.
\end{aligned}
\end{equation*}

The proof concludes by noting that the optimum solutions to the optimization problem \eqref{eq:phi_Jury}
are attained at the boundary of its feasible set. Therefore, having a zero slack variable, \ie, $s_i^\star = 0$, is a necessary condition for $\vert \phi_{2n-s} \vert = r_i^\star$.
\subsection{Proof of Proposition~\ref{prop:alpha_larger_than_one}}\label{app:prop:alpha_larger_than_one}
Recalling that $\vert \phi_{2n-s} \vert$ is characterized by~\eqref{eq:convergence_factor_M}, the proof follows by showing that the inequalities
\begin{enumerate}[label=\roman*)]
\item $g_r^-(1, \beta, \lambda)<g_r^+(1, \beta, \lambda)$ \label{enum:alpha_1}
\item $g_r^+(1, \beta, \lambda) < \max\{g^+_r(\alpha, \beta,\lambda) ,\; g^-_r(\alpha, \beta,\lambda) \}$ \label{enum:alpha_2}
\item $g_1(1, \beta) < g_1(\alpha, \beta)$\label{enum:alpha_3}
\item $g_c(1, \beta, \lambda) < g_c(\alpha, \beta, \lambda)$\label{enum:alpha_4}
\end{enumerate}
hold for $\alpha<1$, $\beta\in(0,\,1)$, and $\lambda\in \{\lambda_i\}_{i\leq n-s}$.

The first inequality~\ref{enum:alpha_1} can be rewritten as $-\beta\lambda < 1$, which holds since $\lambda\geq-1$. As for the second inequality~\ref{enum:alpha_2}, it suffices to show $\Delta g^+ \triangleq g_r^+(1, \beta, \lambda) - g^+_r(\alpha, \beta,\lambda)<0$. After some derivations, we obtain
\begin{equation*}
\begin{aligned}
\Delta g^+=\dfrac{1-\alpha}{2}\left( \sqrt{\lambda^2\beta^2 - 2\beta + 1 + s^+_r} + \lambda \beta -1 \right)
\end{aligned}
\end{equation*}
and observe that $\Delta g^+ < 0$ holds if $\sqrt{\lambda^2\beta^2 - 2\beta + 1 + s^+_r} + \lambda \beta -1 < 0$.
The latter inequality holds for $\lambda \in[-1, 1)$, hence we conclude that $g^+_r(1, \beta,\lambda) < g^+_r(\alpha, \beta,\lambda) \leq \max\{g^+_r(\alpha, \beta,\lambda) ,\; g^-_r(\alpha, \beta,\lambda) \}$.

Next we consider the third inequality~\ref{enum:alpha_3}. For $\alpha\leq 1$ we have $g_1(\alpha, \beta) = 1-\alpha(1-\beta)$. It directly follows that $g_1(1, \beta) < g_1(\alpha, \beta)$ for $\alpha<1$, since $g_1(1, \beta) - g_1(\alpha, \beta) = \alpha - 1$.

In the last step of the proof we address~\ref{enum:alpha_4}. In particular, since $g_c(\alpha, \beta, \lambda)$ is positive, having $\Delta g_c \triangleq g_c(1, \beta, \lambda)^2 - g_c(\alpha, \beta, \lambda)^2 <0$ is equivalent to~\ref{enum:alpha_4}. Thus we study the sign of $\Delta g_c =(1-\alpha)\left(\dfrac{1}{2}\beta\alpha(1-\lambda)  + \dfrac{1}{2}\lambda\beta+ \dfrac{1}{2}\beta - 1\right) + s_c(1) - s_c(\alpha)$.
Using the equality $\dfrac{1}{2}\beta(1-\lambda)  + \dfrac{1}{2}\lambda\beta+ \dfrac{1}{2}\beta - 1= \beta-1$ and $1-\lambda > 0$, for $\alpha <1$ we have
\begin{equation*}
\begin{aligned}
\Delta g_c <& (1-\alpha)\left(\beta-1 \right)+ s_c(1) - s_c(\alpha).
\end{aligned}
\end{equation*}
Recall from Theorem~\ref{thm:reduced_max} that we can only have $\vert \phi_{2n-s}(\alpha, \beta, \lambda) \vert = g_c(\alpha, \beta, \lambda)$ when $s_c(\alpha)=0$. Note that the case when $s_c(\alpha)>0$  corresponds to
\[\vert \phi_{2n-s}(\alpha, \beta, \lambda) \vert = \max\{g^+_r(\alpha, \beta,\lambda) ,\; g^-_r(\alpha, \beta,\lambda),\; g_1(\alpha, \beta) \},\]
which is covered in the previous part of the proof. In the following we let $s_c(\alpha)=0$ and derive the upper bound $s_c(1) < -(1-\alpha)(\beta-1)$. Given the definition of $s_c(1) = \max\{0, -(1/2\beta(1-\lambda) + \beta\lambda )\}$ in Theorem~\ref{thm:reduced_max}, the latter upper bound holds if the following inequalities are satisfied: $(1-\alpha)(\beta-1) <0$ and $\Delta s_c\triangleq(1-\alpha)(\beta - 1) - ( \beta(1-\lambda)/2 + \beta\lambda ) < 0$. The proof concludes by observing that, for $\alpha<1$, $\beta\in(0,1)$, and $\lambda\in[-1,\; 1)$, the former inequality holds, which in turn satisfies the latter inequality, since $\Delta s_c < - ( \beta(1-\lambda)/2 + \beta\lambda )= - 1/2\beta(  1 + \lambda ) < 0$.

\subsection{Proof of Theorem~\ref{thm:optimal_f_alpha}}\label{app:thm:optimal_f_alpha}
Some preliminary results are derived before proving the theorem.
\begin{lem}
\label{lem:consensus_monotonicity_f_complex}
For a fixed $\alpha\in[1,\,2]$, $\lambda\in \{\lambda_i\}_{i\leq n-s}$, and $s_c = 0$, the function $g_c(\alpha, \beta,\lambda)$, defined in \eqref{eq:function_g_lambda}, is monotonically increasing with $\beta\in(0,\,1)$.
\end{lem}
\begin{IEEEproof}
The derivative of $g_c(\alpha, \beta,\lambda)$ with respect to $\beta$ is
\begin{equation*}
\nabla_{\beta} g_c = \frac{1}{2}(\frac{1}{2} \alpha^2\beta(1-\lambda) +1 -\alpha + \alpha \beta\lambda)^{-1/2}  (\frac{1}{2} \alpha (1-\lambda)   + \lambda),
\end{equation*}
which is nonnegative if and only if $\frac{1}{2} \alpha (1-\lambda)   + \lambda \geq 0$. The inequality can be rewritten as $\alpha \geq \dfrac{-2\lambda}{1-\lambda}$, which holds for all $\alpha\in[1,\,2]$ and $\lambda\geq -1$.
\end{IEEEproof}

\begin{lem}
\label{lem:consensus_monotonicity_rho_real}
For a fixed $\alpha\in(0,2]$, $\lambda\in \{\lambda_i\}_{i\leq n-s}$, and $s^{+}_r = s^{-}_r=0$, the functions $g^+_r(\alpha, \beta,\lambda)$ and $g^-_r(\alpha, \beta,\lambda)$ are monotonically decreasing with respect to $\beta$.
\end{lem}
\begin{IEEEproof}
Considering first $g^+_r(\alpha, \beta,\lambda)$, its derivative with respect to $\beta$ is
\begin{align*}
\nabla_{\beta} g^+_r(\alpha, \beta,\lambda) = \dfrac{\alpha}{2}\left( (\lambda^2\beta -1)(\lambda^2 \beta^2 - 2\beta + 1)^{-1/2}  + \lambda \right).
\end{align*}
Since $\beta \in (0,\,1)$ and recalling from Lemma~\ref{lem:bar_lambda} that $|\lambda|\leq 1$, we have $\lambda^2\beta -1 < 0$ and thus
\begin{align*}
\nabla_{\beta} g^+_r(\alpha, \beta,\lambda) &\leq \dfrac{\alpha}{2}\left( (\lambda^2\beta -1)(\beta^2 - 2\beta + 1)^{-1/2}  + \lambda \right)\\
&= \dfrac{\alpha(\lambda - 1)(1+ \lambda \beta)}{2(1-\beta)} < 0.
\end{align*}
Considering $g^-_r(\alpha, \beta,\lambda)$, we have
$\nabla_{\beta} g^-_r(\alpha, \beta,\lambda) = \dfrac{\alpha}{2}\left( (\lambda^2\beta -1)(\lambda^2 \beta^2 - 2\beta + 1)^{-1/2}  - \lambda \right)$.
Similarly as before, $\nabla_{\beta} g^-_r(\alpha, \beta,\lambda)$ can be upper-bounded by
\begin{align*}
\nabla_{\beta} g^-_r(\alpha, \beta,\lambda) &\leq \dfrac{\alpha}{2}\left( (\lambda^2\beta -1)(\beta^2 - 2\beta + 1)^{-1/2}  - \lambda \right)\\
&= \dfrac{\alpha(\lambda +1 )(\lambda \beta -1)}{2(1-\beta)}\leq 0.
\end{align*}
\end{IEEEproof}

\paragraph*{[Proof of Theorem~\ref{thm:optimal_f_alpha}]}
Recall from Proposition~\ref{prop:alpha_larger_than_one} that the minimizing relaxation parameter $\alpha^\star$ lies in the interval $[1,\, 2]$.

First, suppose that $s^{+}_r=s_c = 0$ and observe that $g^+_r \geq g^-_r$ is equivalent to
\begin{equation}\label{eq:grp_greater_grm}
%
\alpha \leq \dfrac{4}{\eta}
\end{equation}
with
\[\eta=2 -(\lambda_{n-s}+ \lambda_1)\beta + \sqrt{\lambda_1^2\beta^2 - 2\beta + 1 } - \sqrt{\lambda_{n-s}^2\beta^2 - 2\beta + 1 }.\]
Recall that $g_1 = \max\{1-\alpha(1-\beta),\; -1 + \alpha(1-\beta) \}$. For $\beta\leq 1/2$, $g^-_r\geq -1 +\alpha(1-\beta)$. Therefore, given \eqref{eq:grp_greater_grm}:
\begin{equation}
|\phi_{2n-s}|=\max \{ g^+_r,g_c,1 -\alpha(1-\beta)\}
\end{equation}
Note that $g^+_r$ is decreasing with respect to $\beta$ while $g_c$ and $1 -\alpha(1-\beta)$ are increasing with respect to $\beta$. Hence, $\beta^\star$ satisfies $  g^+_r=\max\{g_c,1 -\alpha(1-\beta)\}$.  Next, we consider the following cases.

\paragraph*{\ref{enum:case1} and \ref{enum:case2}}

Suppose that $\lambda_{n-s}\neq 0$ and $\beta^\star$ is the solution to $g_r^+ = g_c$, yielding $\beta^\star = ({1-\sqrt{1-\lambda_{n-s}^2}})/{\lambda_{n-s}^2}$. We show $\beta^\star$ is the minimizer by deriving $g_r^+(\alpha,\beta^\star,\lambda_{n-s}) \geq 1-\alpha(1-\beta^\star)$. Since $\beta^{\star2}\lambda_{n-s}^2 - 2\beta^\star + 1 = 0$, the latter inequality can be rewritten as $\beta^\star \leq 1/(2-\lambda_{n-s})$, which is equivalent to $1-\lambda_{n-s}^2 \geq \left(1-\frac{\lambda_{n-s}^2}{2-\lambda_{n-s}} \right)^2$. After manipulations, the former condition reduces to $\lambda_{n-s}(1-\lambda_{n-s})\geq 0$, which holds for $1>\lambda_{n-s} \geq 0$. Hence the minimizing $\beta^\star$ occurs for $g_r^+ = g_c$. Moreover, note that $\beta^{\star2}\lambda_{n-s}^2 - 2\beta^\star + 1 = 0$, which ensures $s^{+}_r = s_c = 0$.

We now fix $\beta^\star$ and optimize over the relaxation parameter. Observing that $g_r^+$ is decreasing with $\alpha$, since $\nabla_\alpha g_r^+(\alpha, \beta^\star, \lambda_{n-s}) = 1/2(-1+\beta^\star\lambda_{n-s})<0$, we conclude that $\alpha^\star$ is the upperbound  in \eqref{eq:grp_greater_grm}.

For the case where $\lambda_{n-s}\geq |\lambda_1|$ and $\lambda_{n-s}\neq 0$ (\ref{enum:case1}), since $\lambda_{n-s} + \lambda_1 \geq 0$ for any choice of $\beta$, $\sqrt{\lambda_1^2\beta^2 - 2\beta + 1} - \sqrt{\lambda_{n-s}^2\beta^2 - 2\beta + 1} \leq 0$,  the upperbound of \eqref{eq:grp_greater_grm} will be larger than 2. On the other hand, $\alpha \in (0,2]$. Thus, $\alpha^\star=2$.

For the case where $ |\lambda_1| \geq \lambda_{n-s}>0$ (\ref{enum:case2}), following a similar line of reasoning to the previous case, we obtain

$$\alpha^\star=\dfrac{4}{2-(\lambda_{n-s} + \lambda_1)\beta^\star + \sqrt{\lambda_1^2 \beta^{\star 2} - 2\beta^\star +1} }\leq 2.$$

\paragraph*{\ref{enum:case3}} As before $\vert\phi_{2n-s}\vert = \max\{g_r^+, g_c, 1-\alpha(1-\beta)\}$. Given lemmas~\ref{lem:consensus_monotonicity_f_complex} and \ref{lem:consensus_monotonicity_rho_real}, the minimizer $\beta^\star$ occurs for $g_r^+ = \max\{g_c, 1-\alpha(1-\beta)\}$. Supposing that the minimizer occurs for $g_r^+(\alpha, \beta^\star, \lambda_{n-s}) =  1-\alpha(1-\beta^\star)$, we obtain $\beta^\star = 1/2$ and $g_r^+(\alpha, \beta^\star, \lambda_{n-s}) =  1-\alpha/2$. Next we show that $1-\alpha/2 \geq g_c(\alpha,\beta^\star,\lambda_{n-s})$, which can be rewritten as
 $ 0 \geq  \frac{\alpha}{2} \lambda_{n-s}  (1-  \frac{\alpha}{2}) + s_c$. Recalling from Theorem~\ref{thm:reduced_max} that $s_c = 0$ is a necessary condition for $\vert\phi_{n-s}\vert = g_c$, we set $s_c$ to zero. Since $\lambda_{n-s}\leq 0$, we have that $g_r^+(\alpha, \beta^\star,\lambda_{n-s}) \geq g_c(\alpha,\beta^\star,\lambda_{n-s})$ for $s_c=0$.

 Next we fix $\beta^\star = 1/2$ and optimize over $\alpha\in[1,\, 2 ]$. Note that $g_r^+(\alpha, \beta^\star, \lambda_{n-s}) =  1-\alpha/2$ is decreasing with $\alpha$ and recall that $\alpha$ is constrained to satisfy~\eqref{eq:grp_greater_grm}. Hence the best parameter $\alpha^\star$ occurs at the boundary of~\eqref{eq:grp_greater_grm}, \ie,  $\alpha^\star = {4}/{(2-\lambda_1)}$, thus concluding the proof.


\subsection{Proof of Lemma~\ref{lem:EWE_Qkappa_general}}\label{app:lem:EWE_Qkappa_general}
Without loss of generality, consider the optimization problem~\eqref{eq:QP_problem} with $\bar{h}=0$ (see Lemma~\ref{lem:QP_change_of_variable}) and include the additional constraint
$\bar{F}^\top( \bar{E} x + \bar{F}) z = 0$.
This constraint may be rewritten as $ z = -\bar{F}^\dagger\bar{E} x $.
Replacing the latter expression in the constraint $\bar{E}x + \bar{F}z = 0$ we obtain $\Pi_{\mathcal{N}(F^\top)}\bar{E}x  = 0$, which can be rewritten as $x=P y$ for some $y\in\mathbf{R}^{s}$. Hence the optimization problem \eqref{eq:QP_problem} is equivalent to
\begin{equation*}
\begin{aligned}
	\underset{y\in\mathbf{R}^{s}}{\mbox{minimize}} &\quad \dfrac{1}{2} y^\top P^\top Q P y + q^\top P y - c^\top \bar{F}^\dagger \bar{E}P y
\end{aligned}
\end{equation*}
The proof follows directly by noting that $P^\top E^\top W E P = P^\top Q P$ yields the same optimal solution of the equivalent problem when $Q\succeq0$ is replaced with $E^\top W E\succ 0$.
\subsection{Proof of Lemma~\ref{lem:feasible}}\label{app:lem:feasible}
Recall that Assumption~\ref{assum:feasible} states that $R$ is chosen so that all solutions to $R(Ex+Fz) = 0$ satisfy $Ex+Fz = 0$. Decomposing $Ex$ as $Ex = \Pi_{\Null{F^\top}}Ex + \Pi_{\mathcal{R}({F})}Ex$, the first equation becomes $R\left(\Pi_{\Null{F^\top}}Ex + \Pi_{\mathcal{R}({F})}(Ex + Fz) \right) = 0$. Since $\Null{F^\top}$ and $\mathcal{R}({F})$ are orthogonal complements, the latter equation can be rewritten as
\begin{subequations}\label{eq:feasible}
\begin{align}
0 & =
RF\left((F^\top F)^{-1}F^\top Ex + z \right)\label{eq:feasible1}\\
0 & = R\Pi_{\Null{F^\top}}Ex.\label{eq:feasible2}
\end{align}
\end{subequations}
The equation~\eqref{eq:feasible1} admits the same solutions as its unscaled counterpart with $R=I$ if and only if $RF$ has an empty null-space, which is equivalent to have $F^\top W F \succ 0$.

As for equation~\eqref{eq:feasible2}, assuming the latter inequality holds and decomposing $REx=\bar{E}x$ as $\bar{E}x = \Pi_{\Null{\bar{F}^\top}}\bar{E}x + \Pi_{\mathcal{R}(\bar{F})}\bar{E}x$, the scaled equations~\eqref{eq:feasible} can be rewritten as
\begin{subequations}
\begin{align}
0 & = \Pi_{\mathcal{R}(\bar{F})}\bar{E}x+\bar{F}z \label{eq:feasible_scaled1}\\
0 & = \Pi_{\Null{\bar{F}^\top}}\bar{E}x,\label{eq:feasible_scaled2}
\end{align}
\end{subequations}
Solutions to~\eqref{eq:feasible2} with $R=I$ can be parameterized as $x=P w\in\Null{\Pi_{\mathcal{N}(F^\top)}E}$, where $P\in\mathbf{R}^{n\times s}$ is an orthonormal basis for $\Null{\Pi_{\mathcal{N}(F^\top)}E}$. Moreover, note that $x=Pw$ is also a solution to~\eqref{eq:feasible_scaled2}, yielding $\Pi_{\Null{\bar{F}^\top}}\bar{E}Pw = 0$. Decomposing $x$ as $x=Pw + P_1 y$,~\eqref{eq:feasible_scaled2} becomes $\Pi_{\Null{\bar{F}^\top}}\bar{E}(Pw + P_1 y) = \Pi_{\Null{\bar{F}^\top}}\bar{E}P_1 y=0$. Thus, \eqref{eq:feasible2} and~\eqref{eq:feasible_scaled2} admit the same solutions $x=Pw$ if and only if
$P_1^\top\bar{E}^\top\Pi_{\Null{\bar{F}^\top}}\bar{E}P_1 \succ 0$.

The proof concludes by observing that, under Assumption~\ref{assum:feasible},~\eqref{eq:feasible2} with $R=I$ and~\eqref{eq:feasible_scaled2} admit the same solutions.


\subsection{Proof of Lemma~\ref{lem:LMI}}\label{app:lem:LMI}
First, suppose that $W=R^\top R$ is chosen such that Assumption~\ref{assum:feasible} holds, as per Lemma~\ref{lem:feasible}. Therefore, we have $\Null{\Pi_{\mathcal{N}(F^\top)}E}$ = $\Null{\Pi_{\mathcal{N}(\bar{F}^\top)}\bar{E}}$.
Note that the unit generalized eigenspace of $(\bar{E}^\top\left(2\Pi_{\mathcal{R}(\bar{F})} - I \right)\bar{E} , \bar{E}^\top \bar{E})$ is characterized by the solutions of the equation $(\bar{E}^\top\left(2\Pi_{\mathcal{R}(\bar{F})} - I \right)\bar{E} - \bar{E}^\top \bar{E})v = 0 $ and corresponds to $\Null{ \Pi_{\mathcal{N}(\bar{F}^\top)}\bar{E}}$.
Hence, we have $P_1^\top\left( \bar{E}^\top\left(2\Pi_{\mathcal{R}(\bar{F})} - I \right)\bar{E} \right)P_1 \prec \lambda_{n-s} P_1^\top (\bar{E}^\top \bar{E})P_1$
and conclude that $\lambda > \lambda_{n-s}$ holds if and only if
\begin{equation}\label{eq:orth_const_P}
P_1^\top\left( \bar{E}^\top\left( 2\Pi_{\mathcal{R}(\bar{F})} - I \right)\bar{E} - \lambda \bar{E}^\top \bar{E} \right)P_1 \prec 0.
\end{equation}
Using the Schur lemma, \eqref{eq:orth_const_P} can be rewritten as~\eqref{eq:upperbound_lambda}.

To conclude the proof, we show that a feasible $W$ with $\lambda \leq 1$ does indeed satisfy the conditions in Lemma~\ref{lem:feasible}. Suppose that~\eqref{eq:orth_const_P} holds with some $\lambda \leq 1$. The inequality $F^\top W F\succ 0$ is clearly satisfied. As for the condition $P_1^\top \bar{E}^\top\Pi_{\Null{\bar{F}}}\bar{E} P_1\succ0$, note that $W$ satisfies
 \iftoggle{draft}{%
\[P_1^\top\left( \bar{E}^\top\left( 2\Pi_{\mathcal{R}(\bar{F})} - I \right)\bar{E} - \lambda \bar{E}^\top \bar{E} - (1-\lambda)\bar{E}^\top \bar{E} \right)P_1 \prec 0,\]
}{%
$P_1^\top\left( \bar{E}^\top\left( 2\Pi_{\mathcal{R}(\bar{F})} - I \right)\bar{E} - \lambda \bar{E}^\top \bar{E} - (1-\lambda)\bar{E}^\top \bar{E} \right)P_1 \prec 0$,
}
since $(1-\lambda)\bar{E}^\top \bar{E} \succeq 0$ for $\lambda\leq 1$. Observing that the latter condition can be rewritten as $2 P_1^\top \bar{E}^\top(I - \Pi_{\Range{F}})\bar{E} P_1=2P_1^\top \bar{E}^\top\Pi_{\Null{\bar{F}}}\bar{E} P_1\succ0$ concludes the proof.

\end{document}